\documentclass[a4paper]{article}

\usepackage{graphbox}
\usepackage{amsfonts}
\usepackage{graphicx,epstopdf}
\usepackage[most]{tcolorbox}

\newcommand{\A}{\mathbf A}
\newcommand{\bb}{\mathbf b}
\newcommand{\M}{\mathbf M}
\newcommand{\nn}{\mathbf n}
\newcommand{\T}{\mathbf T}
\newcommand{\Talg}{\mathbf T_{alg}}
\newcommand{\Tss}{\mathbf{T}_{ss}}
\newcommand{\tcd}{T_f^{c,d}}
\newcommand{\tpq}{t_{pq}}
\newcommand{\spq}{s_{pq}}
\newcommand{\pp}{\mathbf{e}_1}
\newcommand{\qq}{\mathbf{e}_2}
\newcommand{\Ts}{\T_s}
\newcommand{\Xt}{\mathbf X_t}
\newcommand{\X}{\mathbf X}
\newcommand{\Xs}{\mathbf X_s}
\newcommand{\Xss}{\mathbf X_{ss}}
\newcommand{\Tt}{\mathbf T_t}

\begin{document}

\title{On the Evolution of the Vortex Filament Equation for regular $M$-polygons with nonzero torsion}

\author{Francisco de la Hoz$^1$ \and Sandeep Kumar$^2$ \and Luis Vega$^{2,3}$}

\date{
	\footnotesize
	$^1$Department of Applied Mathematics and Statistics and Operations Research, Faculty of Science and Technology, University of the Basque Country UPV/EHU, Barrio Sarriena S/N, 48940 Leioa, Spain. \\[1em]
	$^2$BCAM - Basque Center for Applied Mathematics, Alameda Mazarredo 14, 48009 Bilbao, Spain.\\[1em]
	$^3$Department of Mathematics, Faculty of Science and Technology, University of the Basque Country UPV/EHU, Barrio Sarriena S/N, 48940 Leioa, Spain.
}

\maketitle

\begin{abstract}
In this paper, we consider the evolution of the Vortex Filament equation (VFE):
\begin{equation*}
\X_t = \Xs \wedge \Xss,
\end{equation*}

\noindent taking $M$-sided regular polygons with nonzero torsion as initial data. Using algebraic techniques, backed by numerical simulations, we show that the solutions are polygons at rational times, as in the zero-torsion case. However, unlike in that case, the evolution is not periodic in time; moreover, the multifractal trajectory of the point $\X(0,t)$ is not planar, and appears to be a helix for large times.

These new solutions of VFE can be used to illustrate numerically that the smooth solutions of VFE given by helices and straight lines share the same instability as the one already established for circles. This is accomplished by showing the existence of variants of the so-called Riemann's non-differentiable function that are as close to smooth curves as desired, when measured in the right topology. This topology is motivated by some recent results on the well-posedness of VFE, which prove that the selfsimilar solutions of VFE have finite renormalized energy.

\end{abstract}

\medskip

\noindent\textit{Keywords:}

\noindent Multifractality, Talbot Effect, Turbulence, Vortex Filament Equation

\section{Introduction}

\label{sec:Intro}
Given an initial curve $\X_0 : \mathbb{R} \rightarrow \mathbb{R}^3$, we consider the binormal flow 
\begin{equation}\label{eq:Binorm-flow}
\Xt = \kappa\, \bb,
\end{equation}
where $\kappa$ is the curvature and $\bb$ is the binormal component of the Frenet-Serret system
\begin{equation}\label{mat:SF-frame}
\begin{pmatrix} \T \\ \nn  \\ \bb  \end{pmatrix}_s = 
\begin{pmatrix} 0 & \kappa & 0\\ \mp \kappa & 0 & \tau \\ 0 & -\tau & 0 \end{pmatrix} 
\begin{pmatrix} \T \\ \nn \\ \bb \end{pmatrix}. 
\end{equation}

\noindent The binormal flow can be expressed as 
\begin{equation}\label{eq:VFE}
\X_t = \Xs \wedge_{\pm} \Xss,
\end{equation}

\noindent where $s$ is the arc-length parameter, and $\wedge_\pm$ is the usual cross product with positive and negative sign, corresponding respectively to the Euclidean and hyperbolic geometries. In this paper, we will consider only the positive sign, and, for simplicity's sake, we will omit writing it. (\ref{eq:Binorm-flow}) and its equivalent version (\ref{eq:VFE}) are also referred to as the vortex filament equation (VFE), denomination which we will use along this paper, or as the local induction approximation (LIA). The latter name comes from the local approximation of the Biot-Savart integral that is used to obtain \eqref{eq:Binorm-flow} from the Euler equations (see \cite[p. 510]{Batchelor} and \cite{saffman}).

We can also consider the tangent vector $\T=\Xs$. Then, differentiating \eqref{eq:VFE} with respect to $s$ yields 
\begin{equation}\label{eq:SMP}
\Tt = \T \wedge \Tss,
\end{equation}

\noindent which is the Schr\"{o}dinger map equation onto the sphere $\mathbb S^2$. Moreover, it follows immediately that $\T$ preserves its length, so we can assume after a scaling that $\T \in \mathbb{S}^2$.  In \cite{hasimoto}, Hasimoto made a fundamental contribution by defining the filament function:
\begin{equation}\label{eq:Hasimoto}
\psi(s,t) = \kappa(s,t) e^{i\int_0^s \tau(s^\prime,t) ds^\prime}, 
\end{equation}

\noindent and proving that such a $\psi$ solves the nonlinear Schr\"{o}dinger (NLS) equation:
\begin{equation}\label{eq:NLS}
\psi_t = i \psi_{ss} + i \left(\frac{1}{2}(|\psi|^2 + A(t)) \right)\psi, \quad A(t) \in \mathbb{R}.
\end{equation}

\noindent Coming back to VFE, some of its explicit solutions are the circle, the helix, and the straight line. The last case corresponds to the line vortex, which does not move \cite[p. 93]{Batchelor}. As for the circle, it moves normal to its plane with a constant velocity given by the inverse of the radius, and the choice of the direction is determined by the orientation given to the circle; this example represents the smoke rings that are generated in real fluids \cite[p. 522]{Batchelor}. Finally, the helix does not change its shape either, but it screws up or down, depending also on the orientation. Vortices with a helical shape are easy to generate; for example, they can be seen behind the tips of the blades of a propeller. In fact, there are explicit solutions of Euler equations that have a helical shape \cite{Hardin,Ricca}.

In \cite{HozVega2014}, the case of regular planar polygons of $M$ sides with zero torsion (from now on, referred to as planar $M$-polygons) was considered. Because of the mix of lack of regularity (due to the presence of corners) and periodicity, the available analytical techniques are not sufficient to treat this type of curves. Recently, it has been proved that, for polygonal lines that, at any given time, are asymptotically close at infinity to two straight lines, the initial value problem is well posed in an appropriate topology \cite{BanicaVega2018}. This asymptotic behavior allows the waves generated by the corners to escape to infinity, so that their interaction is not strong enough to generate new corners. As observed in \cite{Didier2}, this is not the case when closed polygons are considered. As proved in \cite{HozVega2014} for planar $M$-polygons, these waves can interact in a rather complicated way, so new corners appear at rational multiples of the time period $2\pi/M^2$, in a similar way as that exhibited by optical waves in the so-called Talbot effect \cite {Berry,ET,Olver}. As explained in \cite{HozVega2014}, similar effects have been observed in the evolution of noncircular jets \cite{GG}. In particular, if nozzles with the shape of equilateral triangles or squares are considered, some complex structures in the shape of skew polygons are observed at later times (see \cite[Fig. 6]{GG}, and also \cite[Fig. 10 and p. 1492]{GGP}, where it is mentioned ``[...] a consistent eightfold distribution pattern is also suggested [...]'', for square nozzles). 
 
This complicated dynamics is also present when one looks at the trajectory of any of the corners of a planar $M$-polygon. The trajectory seems to be a multifractal, very reminiscent of the graph of Riemann's non-differentiable function:
 \begin{equation} \label{eq:RND}
\sum_{k=1}^{\infty} \frac{\sin{\pi  k^2 t}}{ \pi k^2}, \quad t \in [0,2].
\end{equation}

\noindent Indeed, in \cite[Section 5.1]{HozVega2014}, it was observed that, as the polygonal curve $\X(s,t)$ evolves in time, it also moves in the vertical direction with constant speed $c_M$. Additionally, taking into account the symmetries of the problem, it is easy to see that, for a given $M$, the trajectory of a single point $\X(0,t)$ lies in a plane. Hence, removing the vertical movement, and using complex coordinates, we can define
\begin{equation*}
\tilde z_M(t) =  -\| (X_1(0,t), X_2(0,t) ) \|_2 + i (X_3(0,t) - c_M t), \quad t \in [0,2\pi/M^2],
\end{equation*}

\noindent which is a closed, $2\pi/M^2$-periodic curve. Then, using an appropriate scaling that depends on $M$, very strong numerical evidence is given to show that, as $M\to\infty$, $\tilde z_M(t)$ converges to a complex version of Riemann's non-differentiable function, namely
\begin{equation} \label{eq:phi}
\phi(t) = \sum_{k=1}^{\infty} \frac{e^{\pi i k^2 t}}{i \pi k^2}, \quad t \in [0,2].
\end{equation} 

\noindent We refer the reader to \cite{HozVega2014} for the details, to \cite{Ja} for the proof of the multifractality of Riemann's function, and to \cite{CC}, where other relevant trigonometric sums are considered. On the other hand, since the planar $M$-polygon tends initially to a circle as $M\to\infty$, it is proved as a consequence that the time evolution of a circle is not stable, at least at the numerical level. In other words, a particle can be located in a curve arbitrarily close to a circle, but in the right topology, its trajectory converges to the multifractal graph given by $\phi$.

Later on, in \cite[Section 6]{HozVega2018}, the linear momentum of $\X$ (sometimes also called impulse in the fluid literature) was studied. In particular, it was observed that it is a periodic function with a period depending on $M$. After renormalizing the period, a spectral analysis was done and it was proved numerically that the leading Fourier coefficients behave as those in \eqref{eq:RND}. A similar analysis can be done for $\tilde z_M(t)$, obtaining
\begin{equation*}
 \lim\limits_{M\rightarrow \infty} |n\,\tilde a_{n,M}|  = 
 \begin{cases}
 1, & \text{if } n=k^2,\ k\in\mathbb N, \\
 0, & \text{otherwise};
 \end{cases} 
 \end{equation*}

\noindent where $\tilde a_{n,M}$ are the Fourier coefficients of the scaled $\tilde z_M(t)$, and $\mathbb N = \{1, 2, 3, \ldots\}$.
 
The main purpose of this paper (described in Section \ref{sec:X0t}) is to prove that helices and straight lines also have the same kind of instability as the circle. We will show this by approximating them using regular $M$-polygons with nonzero torsion (from now on, referred to as helical $M$-polygons), i.e., regular polygons whose tangent vectors are $2\pi$-periodic functions that take just $M$ values on $\mathbb S^2$. As a consequence, we can say that Riemann's non-differentiable function and its variants appear as universal objects in the dynamics of singular solutions of \eqref{eq:VFE}. Moreover, these universal objects are, in the right topology, as close as desired to smooth curves.
 
The paper is organized as follows. After defining the problem in Section \ref{sec:prob-def}, we introduce in Section \ref{sec:prob-form} the parametric form of our initial data, i.e., the curve $\X$ and its tangent vector $\T$, whose third component is a constant denoted by $b\in(0,1)$. Note that $b=0$ reduces back to the planar $M$-polygon case, $b=1$ to the straight line, and, for an intermediate value of $b$, the corresponding polygonal curve has a helical shape. Thus, the initial data is characterized by two parameters: $M$ (which reflects the periodicity of $\T$) and $b$. These two parameters also determine the curvature angle $\rho_0$ and the torsion angle $\theta_0$. On the other hand, at the level of the NLS equation, the nonzero-torsion problem can be seen as a Galilean transformation of the planar $M$-polygon problem. Thus, denoting \eqref{eq:Hasimoto} as $\psi_\theta$, when $\theta_0\neq0$, and as $\psi$, when $\theta_0=0$, the initial data for the two problems can be related by
\begin{equation*}
\psi_{\theta}(s,0) = C\, e^{i\gamma s} \psi(s,0),
\end{equation*}
where  $\gamma = M\theta_0/2\pi$, and $C$ is a constant depending on the initial structure of the $M$-polygons in both cases. Then, using the Galilean invariance of \eqref{eq:NLS}, we obtain 
\begin{equation}
\label{e:psithetapsi}
\psi_{\theta}(s,t) = C \, e^{i\gamma s- i \gamma^2 t} \psi(s-2\gamma t,t).
\end{equation}

\noindent In Section \ref{sec:psi-comp}, we compute $\psi_\theta(s,t)$, at $t=\tpq=  (2\pi/M^2)(p/q)$, with $\gcd(p,q)=1$, i.e., at times that are rational multiple of $2\pi/M^2$. Since $\psi(s,\tpq)$ is known from the zero-torsion case, it follows from \eqref{e:psithetapsi} that the helical $M$-polygon curve at time $\tpq$ is a polygon with $Mq$ sides (if $q$ odd) or $Mq/2$ sides (if $q$ even). The structure of the new polygon can be determined by the generalized quadratic Gau{\ss} sums $G(-p,m,q)$, as in \cite{HozVega2014}; but the presence of torsion causes a lack of space and time periodicity. In other words, unlike in the planar polygon case, we observe, at the end of one time-period, a polygonal curve with the same number of sides in a spatial period, but with the corners moved by a certain amount, which we call Galilean shift, and rotated by a certain angle, which we call phase shift. Both shifts depend on the torsion introduced in the problem. 

In Section \ref{sec:Alg-sol}, we construct the algebraic solution at any rational time $\tpq$. The arguments on the Dirac deltas in (\ref{eq:psi-theta-st_final}) determine the Galilean shift, and thus, we know the exact location of the new corners. Then, a principle of conservation of energy proved in \cite{BanicaVega2018} allows us to compute the angle $\rho_q$ between any two adjacent sides of the resulting polygon, so, using the generalized quadratic Gau{\ss} sums, we can construct the polygonal curve, up to a rigid movement. Later on, using algebraic techniques as in \cite{HozVega2018}, an expression for the speed of the center of mass $c_M$ is given. With all these ingredients, the algebraic solution can be computed up to a rotation that remains undetermined. 

In Section \ref{sec:Num-simu}, we comment on the numerical method to solve (\ref{eq:VFE})-(\ref{eq:SMP}); as in \cite{HozVega2014}, we are able to make use of the symmetries of the tangent vectors and reduce the computation cost quite effectively. As already said, it was noted in \cite{HozVega2014} that, as the polygonal curve evolves in time, it also moves in the vertical direction with constant speed. A similar phenomenon is observed for the helical polygon as well; we compute the angle between any two adjacent sides of the polygonal curve, the speed of the center of mass, and compare those values with their algebraic counterparts. 

In Section \ref{sec:X0t}, we study the trajectories of one corner initially located at $s=0$. The helical case is studied in Section \ref{0<b<1}. During the time evolution, the numerical simulations suggest that, for a given $M$, and $b\in(0,1)$, besides a constant vertical movement, the polygonal curve rotates around the $z$-axis. The trajectory of a single point is a multifractal, but no more planar. Moreover, by taking $b$, such that $\theta_0 = \pi c/d$, with $\gcd(c,d)=1$, $c, d \in \mathbb N$, the periodicity in space can be recovered. Furthermore, for such time period, after removing the vertical height, the third component of $\X(0,t)$ is periodic, and its structure can be compared with the imaginary part of
$$
\phi_{c,d}(t) = \sum_{k \in {A_{c,d}}} \frac{e^ {{2\pi i k t}}}k, \quad t\in
\begin{cases}
	[0,1/2], & \mbox{if $c\cdot d$ odd,} \\
	[0,1], & \mbox{if $c\cdot d$ even,} 
\end{cases}
$$
where the set $A_{c,d}$ is defined in \eqref{eq:dom-pts}. After applying a scaling as before, and expressing it in terms of its Fourier expansion, we get
\begin{equation*}
	\lim\limits_{M\to \infty} |n\, b_{n,M}|  = 
	\begin{cases}
		1/4, & \text{if $n \in A_{c,d}$ and $c\cdot d$ odd},
		\\
		1/2, & \text{if $n \in A_{c,d}$ and $c\cdot d$ even},
		\\
		0, & \text{otherwise}.
	\end{cases} 
\end{equation*}

\noindent where $b_{n,M}$ are the Fourier coefficients. 

In the case of the straight lines, studied in Section \ref{b=1}, $b$ must be taken very close to $1$. Here, the trajectory $\X(0,t)$ in the XY-plane tends to a $2\pi$-periodic closed curve, which can be compared to
$$
\phi_M(t) =  \sum\limits_{k \in A_M} \frac{e^{2\pi i k^2 t}}{k^2}, \quad t \in [0,1],
$$
where the set $A_M$ is defined in \eqref{e:AM}. Thus, we consider first
$$
z(t) = \frac{ X_1(0,t)}{1+ \tilde X_3(t)} + i \frac{X_2(0,t)}{1+\tilde X_3(t)},
$$
with $\tilde X_3(0,t)=X_3(0,t)-c_M t$, $t\in[0,2\pi]$, and after rotating it counterclockwise $\pi / 2- \pi / M$ radians, define $z_M(t)$. As before, we approximate the Fourier coefficients $c_{n}$ of a properly scaled version of $z_M$, to show that, for a given $M$,
\begin{equation*}
\lim\limits_{b\to 1^{-}} |n \,c_{n}|  = 
\begin{cases}
1, & \text{if } n \in A_{M}, \\
0, & \text{otherwise}.
\end{cases} 
\end{equation*}

\noindent In Section \ref{sec:Talg-irrational}, we examine the behavior of $\T$ for rational times $t_{pq}$, with $q\gg1$, and compare it with the zero-torsion case. Finally, Section \ref{sec:M-1-relation} describes briefly the numerical relationship between the $M$-corner and the one-corner problems. Following the same steps as in \cite{HozVega2018}, we try to answer how the helical $M$-polygon problem can be understood as a superposition of $M$ one-corner problems at infinitesimally small times.

\section{A solution of $\Xt = \Xs \wedge \Xss$ for a regular helical $M$-polygon}

The first goal of this paper is to construct the solutions  of (\ref{eq:VFE}) and describe their corresponding dynamics, for initial data given by regular non-planar $M$-polygons. As will be explained in this section, the behavior of regular polygons with nonzero torsion is a consequence of the Galilean symmetry present in the set of solutions of \eqref{eq:NLS}. Therefore, this paper can be regarded an extension of the zero-torsion case considered in \cite{HozVega2014}, and as a result, \cite[Th. 1]{HozVega2014} holds true here as well, only needing to change \cite[(23)]{HozVega2014} and \cite[(25)]{HozVega2014} by \eqref{eq:psi-theta-st_final} and \eqref{eq:rho_q} below, respectively.

\begin{figure}[!htbp]\centering
	\includegraphics[width=0.5\textwidth, clip=true]{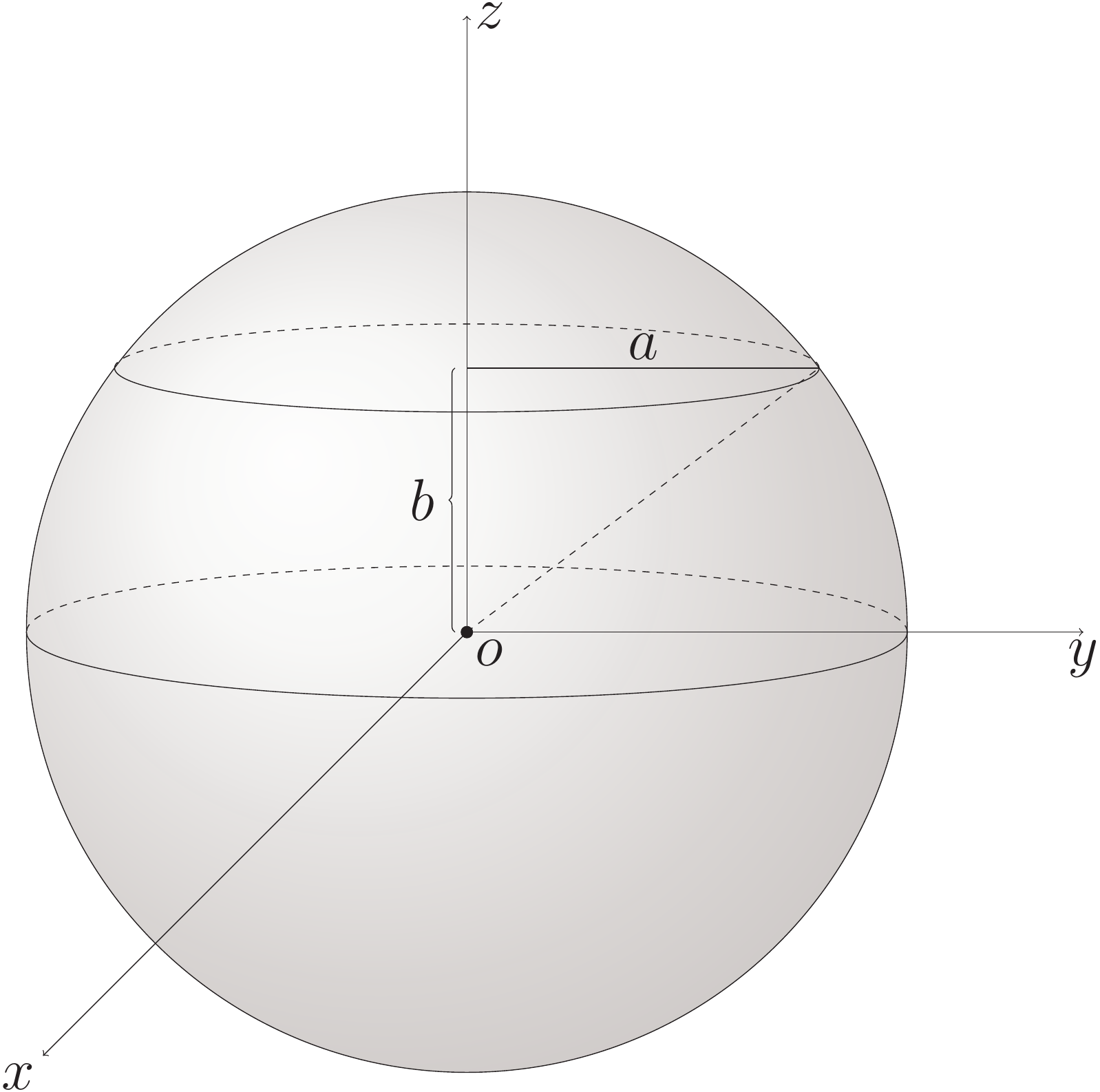}
	\caption{Unit sphere $\mathbb{S}^2$, and parameters $a$ and $b$, where $a^2+b^2=1$.}
	\label{fig:Unit-sphere}
\end{figure}

\subsection{Problem definition}\label{sec:prob-def}

Let us consider an arc-length parameterized regular $M$-polygon with torsion depending on a parameter $b$. Due to the fact that (\ref{eq:VFE}) and  (\ref{eq:SMP}) are rotation invariant, we can assume that the $2\pi$-periodic tangent vector $\T(s,0)$ lies on a circle of radius $a$, with $a^2+b^2=1$ (see Figure \ref{fig:Unit-sphere}):
\begin{equation} \label{eq:Ini-data-T}
\T(s,0) = \left(a\cos\left(\frac{2\pi k}{M}\right),a\sin\left(\frac{2\pi k}{M}\right),b \right) \equiv  (a\,e^{2\pi i k/M},b), \quad s\in(s_k,s_{k+1}),
\end{equation}
where $k=0,1,\ldots,M-1$. The corresponding curve $\X(s,0)$ is a helical polygon with corners located at
\begin{equation}\label{eq:Ini-data-X}
\X(s_k,0) = \left(\frac{a\pi\sin(\pi(2k-1)/M)}{M\sin(\pi/M)}, -\frac{a\pi\cos(\pi(2k-1)/M)}{M\sin(\pi/M)}, b\,s_k\right),
\end{equation}

\noindent so that, for any $s\in(s_k,s_{k+1})$, the corresponding point $\X(s,0)$ lies on the line segment joining $\X(s_k,0)$ and $\X(s_{k+1},0)$. Since $\T \in \mathbb{S}^2$, it follows that $b\in[-1,1]$; in this paper, we work with $b>0$, because the case with $b<0$ can be recovered by the symmetries.

Let us denote $\T_k=\T(s,0)$, for $s\in(s_k,s_{k+1})$. The curvature angle $\rho_0$ does not depend on $k$, and it is defined as the angle between $\T_k$ and $\T_{k+1}$, for all $k$:
\begin{equation} \label{eq:rho_expression}
\rho_0 = 2 \arcsin\left(a\sin\left(\frac{\pi}{M}\right)\right).
\end{equation}

\noindent On the other hand, the torsion angle $\theta_0$ does not depend on $k$, either, and it is defined as the angle between $(\T_{k-1} \wedge \T_{k})$ and $(\T_{k} \wedge \T_{k+1})$, for all $k$:
\begin{equation} \label{eq:theta_expression}
\theta_0 = 2 \arctan\left(b\tan\left(\frac{\pi}{M}\right)\right).
\end{equation}
Observe that both angles satisfy
\begin{align*}
\cos\left(\frac{\rho_0}{2}\right)\cos\left(\frac{\theta_0}{2}\right) & = \cos\left(\arcsin\left(a\sin\left(\frac{\pi}{M}\right)\right)\right)\cos\left(\arctan\left(b\tan\left(\frac{\pi}{M}\right)\right)\right)
	\\
& = \sqrt{\frac{1 - a^2\sin^2\left(\frac{\pi}{M}\right)}{1 + b^2\tan^2\left(\frac{\pi}{M}\right)}} = \cos\left(\frac\pi M\right).
\end{align*}

\subsubsection{Spatial symmetries and the Hasimoto transformation}

\label{sec:Spatial_Sym_X_T}

The symmetries of the initial data for the zero-torsion case (see \cite[(29)]{HozVega2014}) are also valid here, and as a result $\X(s,t)$ and $\T(s,t)$ are invariant under a rotation of angle $2\pi k/M$ around the $z$-axis, with the only exception that the third component $X_3$ has now a translation symmetry, i.e., $X_3(s+2\pi k/M,t) = X_3(s,t) + 2 \pi k\, b /M$, for all $k\in \mathbb{Z}$, $t\geq0$. Moreover, $\T(s,t)$ and $\T(-s,t)$ are symmetric about the XZ-plane, and, consequently, $\X(s,t)$ and $\X(-s,t)$ are symmetric about the $y$-axis, for all $t$. An important fact that will be useful later is that $\X(s+2\pi k,t) - \X(s,t) = 2 \pi b\, k (0,0,1)$, for all $k\in \mathbb{Z}$. 

We denote the filament function  (\ref{eq:Hasimoto}) by $\psi_\theta$, for $\theta_0 \neq 0$, and use $\psi$, when $\theta_0=0$. To avoid issues related to vanishing curvatures, we work with a more general form of the Frenet-Serret formulas, where $\T$, $\pp$, $\qq$ are the orthonormal basis vectors \cite{Bishop} (see \cite[Section 2]{HozVega2014}, for more details). However, in order to obtain $\T$, $\pp$ and $\qq$ at any rational time, we have to transform $\psi_\theta$, by multiplying it by a certain constant. 

\subsection{Problem formulation}\label{sec:prob-form}

Let us write the initial data of the NLS equation corresponding to the planar $M$-polygon problem as

$$
\psi(s,0) = c_0 \sum_{k=-\infty}^{\infty} \delta(s-2\pi k/M), \quad s \in [0,2\pi],
$$

\noindent and, when $\theta_0 \neq 0$, as
$$
\psi_{\theta}(s,0) = c_{\theta,0}  \ e^{i\gamma s} \sum_{k=-\infty}^{\infty} \delta(s-2\pi k/M), \quad s \in [0,2\pi],
$$

\noindent where $\gamma = M \theta_0/ 2\pi$ satisfies $\lim_{M\to \infty} \gamma = b$; and $c_0$ and $c_{\theta,0}>0$ are constants depending on the initial configuration of the respective curve. In particular, taking
\begin{equation}
\label{e:ctheta0}
c_{\theta,0} = \sqrt{-\frac2\pi  \ln  \left(\cos  \left(\frac{\rho_0}{2}\right)\right)},
\end{equation}

\noindent we have
\begin{equation}\label{eq:Ini-data-torsion}
\psi_\theta(s,0) = \frac{c_{\theta,0}}{c_0} \ e^{i \gamma s} \psi(s,0).
\end{equation}
At this point, we recall one important symmetry of the NLS equation called the Galilean invariance, which states that, if $\psi$ is a solution of \eqref{eq:NLS}, then, for the initial datum $\psi_{n}(s,0) = e^{i n\,s} \psi(s,0)$, the corresponding solution is given by $\psi_{n}(s,t) = e^{i n\,s - i n^2t} \psi(s-2n\, t,t)$, with $n, t\in \mathbb{R}$, and satisfies \eqref{eq:NLS}, too. Bearing in mind this fact, the solution corresponding to \eqref{eq:Ini-data-torsion} is given by 
\begin{equation} \label{eq:psi-theta}
\psi_{\theta}(s,t) =\frac{c_{\theta,0}}{c_0} \ e^{i\gamma s- i \gamma^2 t} \psi(s-2\gamma\, t,t).
\end{equation}
Hence, using \cite[(33)]{HozVega2014}, we get
\begin{equation} \label{eq:psi-theta-st}
\psi_\theta(s,t) =   \frac{c_{\theta,0}}{c_0}  \hat{\psi}(0,t) e^{i\gamma s - i \gamma^2 t} \sum\limits_{k=-\infty}^{\infty}  e^{-i(Mk)^2 t+iMk(s-2\gamma t)},
\end{equation} 
where $\hat{\psi}(0,t)$ is a constant depending on time, which, as mentioned in \cite{HozVega2014}, can be assumed to be real, for all $t$. Note that, in (\ref{eq:psi-theta}), $\psi$ is $2\pi/M$-space-periodic and $2\pi/M^2$-time-periodic; but, when $\gamma\in(0,1)$, $\psi_\theta(s,t)$ is not periodic. However, by taking $\gamma$ rational, space periodicity can be recovered at large times; and if $\gamma=1$, $\psi_\theta(s,t)$ is both space and time periodic. During the evolution of an $M$-polygon, this will give rise to a phase shift and a Galilean shift, as explained later on in this paper.

\subsection{Computation of $\psi_\theta(s,t)$ for rational multiples of $t=2\pi /M^2$}

\label{sec:psi-comp}

When $t=\tpq=(2\pi/M^2)(p/q)$, with $p \in \mathbb{Z}$, $q \in \mathbb{N}$, and $\gcd(p,q)=1$, we compute $\psi_\theta(s,t_{pq})$ by substituting $t_{pq}$ and $\gamma$ in (\ref{eq:psi-theta-st}). Then, similarly as in \cite{HozVega2014} for $\psi(s,\tpq)$, we have
\begin{align*}
\psi_\theta(s,t_{pq}) & = \frac{2\pi}{Mq} \frac{c_{\theta,0}}{c_0}   \hat{\psi}(0,\tpq) e^{i(\theta_0^2/(2\pi))(p/q)}  \sum\limits_{k=-\infty}^{\infty}   \sum\limits_{m=0}^{q-1}   
	\\
& \qquad\qquad G(-p,m,q)  e^{ i   (k\theta_0+m\theta_0/q)}   \delta\left(s-\frac{2\theta_0 p}{Mq} - \frac{2\pi k}{M} -\frac{2\pi m}{Mq}\right), 
\end{align*}

\noindent where $G(-p,m,q) = \sum_{c=0}^{q-1} e^{(-2\pi i c^2p + 2\pi i  cm)/q}$ is a generalized quadratic Gau{\ss} sum. Using the properties of these sums (see \cite[Appendix]{HozVega2014}), we get
\begin{align}\label{eq:psi-theta-st_final}
\psi_\theta(s,t_{pq}) =
\begin{cases}
\frac{2\pi}{M \sqrt{q}} \frac{c_{\theta,0}}{c_0} \hat{\psi}(0,\tpq) e^{i(\theta_0^2/(2\pi))(p/q)}  \sum\limits_{k=-\infty}^{\infty}
\sum\limits_{m=0}^{q-1}
	\\
\quad e^{ i   (\xi_m+k\theta_0+m\theta_0/q)} \delta\left(s- \frac{2\theta_0 p}{Mq} - \frac{2\pi k}{M} -\frac{2\pi m}{Mq}\right), &\text{if $q$ odd},
	\\
\frac{2\pi}{M \sqrt{q/2}} \frac{c_{\theta,0}}{c_0} \hat{\psi}(0,\tpq) e^{i(\theta_0^2/(2\pi))(p/(q/2))}  \sum\limits_{k=-\infty}^{\infty} \sum\limits_{m=0}^{q/2-1}
	\\
\quad e^{ i   (\xi_{2m}+k\theta_0+2m\theta_0/q) }\delta\left(s-\frac{2\theta_0 p}{Mq} - \frac{2\pi k}{M} -\frac{4\pi m}{Mq}\right),  &\text{if $\frac q2$ even},
	\\
\frac{2\pi}{M \sqrt{q/2}} \frac{c_{\theta,0}}{c_0} \hat{\psi}(0,\tpq) e^{i(\theta_0^2/(2\pi))(p/(q/2))} \sum\limits_{k=-\infty}^{\infty}\sum\limits_{m=0}^{q/2-1} 
	\\
\quad e^{ i   (\xi_{2m+1}+k\theta_0+(2m+1)\theta_0/q) }\delta\left( s-\frac{2\theta_0 p}{Mq} - \frac{2\pi k}{M} -\frac{2\pi  (2m+1)}{Mq}\right),  &\text{if $\frac q2$ odd},
\end{cases}
\end{align}

\noindent for a certain angle ${\xi_m}$ depending on $m$, $p$ and $q$. Hence, at any rational time $t_{pq}$, the initial $M$ Dirac deltas in $s\in[0,2\pi)$ turn into $Mq$ Dirac deltas (if $q$ odd), or $Mq/2$ Dirac deltas (if $q$ even). Moreover, at those times, the absolute value of their coefficients is constant, and since the Dirac deltas are equally spaced, the sides of the resulting polygon are equally-lengthed. On the other hand, as a result of the Galilean transformation, a corner initially located at $2\pi k/M$, $k\in\mathbb Z$, is translated by $\spq=2\theta_0 p/Mq$ at time $\tpq$; we call this extra movement the Galilean shift. Although, strictly speaking, $\psi_\theta$ (and $\X$ and $\T$) are not time-periodic now, their structure repeats whenever $t$ is increased by $2\pi/M^2$; along this paper, we denote this important quantity $T_f \equiv 2\pi/M^2$, and refer to it, with some abuse of language, as the time period.

\subsubsection{Computation of $\hat{\psi}(0,\tpq)$}

\label{sec:rho_q-comp-alg}

In \cite{GutierrezRivasVega2003}, the following expression was given for the one-corner problem:
\begin{equation} \label{eq:rho_c0}
\cos(\rho_0/2) = e^{-\pi c_0^2/2},
\end{equation}
where $\rho_0$ is the angle formed by the tangent vectors at the corner, when the singularity happens, and $c_0$ is the curvature or coefficient of the Dirac delta in the initial solution of the NLS equation. Later on, in \cite{HozVega2018}, it was shown that, for small times, (\ref{eq:rho_c0}) holds true for planar $M$-polygons as well; in particular, at $t=0$,
$$
c_0 = \sqrt{-\frac2\pi \ln  \left(\cos\left(\frac\pi M\right)\right)} \Longleftrightarrow \rho_0 = \frac{2\pi}{M}.
$$

\noindent On the other hand, in \cite{HozVega2014}, $\hat{\psi}(0,\tpq)$ was guessed from the numerical results, for planar $M$-polygons. This was corroborated in \cite[Section 7]{HozVega2018}, where it was also suggested that $\hat{\psi}(0,\tpq)$ might be obtained from a conservation law. In \cite[Section 4]{BV2019}, it has been recently observed that this conservation law is a consequence of the one established for non-closed polygonal lines in \cite[Th. 14]{BanicaVega2018}. The argument is as follows.

Let $\psi(s,0) = \sum_{k} \alpha_k(0) \delta(s-2\pi k/M)$. Then, in \cite{BanicaVega2018}, assuming some decay conditions of $\alpha_k(0)$, the existence of a unique solution $\psi(s,t)= \sum_{k} \alpha_k(t) e^{it\partial^2_s} \delta(s-2\pi k/M)$ of \eqref{eq:NLS} is proved, with $t\,A(t)=\sum_{k} |\alpha_k(0)|^2$. Moreover, $\sum_{k} |\alpha_k(t)|^2=\sum_{k} |\alpha_k(0)|^2$. Assume now that $\alpha_{k+M} = \alpha_k$, for all $k$. Then, if the solution exists and is unique,  $\alpha_{k+M}(t) = \alpha_k(t)$,  and $\sum_{j=0}^{M-1} |\alpha_j|^2 = Q^2 = \sum_{j=0}^{M-1}|\alpha_j(t)|^2$. In our case, at $t=0$,
$$
|\alpha_k|=c_{\theta,0} \Longrightarrow Q^2= \sum_{k=0}^{M-1}|c_{\theta,0}|^2= M c_{\theta,0}^2= -\frac{2M}\pi \ln \cos\left(\frac{\rho_{0}}2\right).
$$

\noindent Furthermore, at any rational time $t_{pq}$ (taking $q$ odd for now), there are from \eqref{eq:psi-theta-st_final} $Mq$ Dirac deltas with coefficients of equal modulus, which we call $c_{\theta,q}$. Therefore, 
\begin{align*} 
 M c_{\theta,0}^2= Q^2 &= \sum\limits_{k=0}^{Mq-1} |c_{\theta,q}|^2 = Mq \left|\frac{2\pi c_{\theta,0}}{M \sqrt{q} c_0} \hat{\psi}(0,\tpq) \right|^2,
\end{align*}

\noindent which yields $ \hat{\psi}(0,\tpq) = Mc_0 / 2\pi$ and $c_{\theta,q} = c_{\theta,0} / \sqrt{q}$. Moreover, $\lim_{M\to \infty} \hat{\psi}(0,\tpq) = 1/2\sqrt{\pi}$. Note that (\ref{eq:rho_c0}) holds true whenever there is a singularity formation; e.g., in our case, at rational times $\tpq$. The expression for $c_{\theta,q}$ shows that the angle $\rho_q$ between any two tangent vectors remains equal and can be computed by writing
\begin{align*}
\sqrt{-\frac{2}{\pi} \ln \left( \cos \left(\frac{\rho_q}{2}\right)\right)}
= \frac{1}{\sqrt{q}}\sqrt{-\frac{2}{\pi} \ln \left( \cos \left(\frac{\rho_0}{2}\right)\right)} 
\implies \cos \left(\frac{\rho_0}{2}\right) = \cos^q\left(\frac{\rho_q}{2}\right). 
\end{align*}

\noindent After proceeding in the same way for $q$ even, we conclude that
\begin{equation}\label{eq:rho_q}
\cos(\rho_q/2) =
\begin{cases}
 \cos^{1/q}(\rho_0/2), \quad & \mbox{if $q$  odd}, \\
 \cos^{2/q}(\rho_0/2), \quad & \mbox{if $q$  even}, 
\end{cases}
\end{equation}

\noindent which is the same expression as in \cite[(25)]{HozVega2014}, with $\rho_0=2\pi/M$.  

\subsection{Algebraic solution}

\label{sec:Alg-sol}

By using the same technique as in \cite{HozVega2014}, $\X$ and $\T$ can be computed algebraically, up to a vertical movement and a rotation. Their algebraic counterparts, denoted respectively by $\X_{alg}$ and $\T_{alg}$, are obtained as follows. Let us define, for any rational time $\tpq$ with $q$ odd (the case with $q$ even is similar):
\begin{equation}\label{eq:Scaling}
\Psi_\theta  (s,\tpq) = \frac{\rho_q}{c_{\theta,q}} e^{-i(\theta_0^2/(2\pi))(p/q)}\psi_\theta(s,\tpq);
\end{equation}
and, since $ \lim_{q\to \infty} \rho_qe^{-i(\theta_0^2/(2\pi))(p/q)}/c_{\theta,q}<\infty$, $\Psi_\theta$ is well-defined. After writing $\Psi_\theta  (s,\tpq) =  (\alpha + i \beta)(s,\tpq)$, we integrate the Frenet-Serret frame and obtain the $Mq$ basis vectors $\tilde{\T}(s)$, $\tilde{\mathbf{e}}_1(s)$, $\tilde{\mathbf{e}}_2(s)$, at $s=  (2\pi  (k+1) /Mq + \spq)^-$, by the action of $Mq$ rotation matrices denoted as $\mathbf{R}_k$, $k=0,1,\ldots,Mq-1$ (see \cite[Section 3.4]{HozVega2014} for the details). Moreover, it is important to mention that 
\begin{equation} \label{mat:torsion-conservation}
\mathbf{R}_{Mq-1} \cdot \mathbf{R}_{Mq-2}  \ldots  \mathbf{R}_1 \cdot \mathbf{R}_0 = \mathbf{R}^x_{M\theta_0} = \begin{pmatrix} 	1 & 0 & 0 \\ 0&\cos (M\theta_0) & -\sin (M\theta_0)  \\ 0 & \sin (M \theta_0)   &  \cos (M \theta_0)   \end{pmatrix},
\end{equation}
which holds true for any $q$ and $M$, and implies that the following quantity is preserved: 
\begin{equation}\label{eq:torsion-conservation}
\int_{0}^{2\pi} \tau(s^\prime,\tpq) ds^\prime = \int_{0}^{2\pi} \tau(s^\prime,0) ds^\prime = \frac{M}{2\pi} \theta_0 \int_{0}^{2\pi}ds^\prime = M \theta_0.
\end{equation}

\noindent On the other hand, $\tilde{\X}$, i.e., $\X$ up to a rigid movement, can be obtained by
\begin{align*}
\begin{cases}
\tilde{\X}(s_{pq}) = \tilde{\X}(0) + s_{pq} \tilde{\T}(s_{pq}^-), \\
\tilde{\X}(\frac{2\pi  (k+1)}{Mq}+s_{pq}) = 	\tilde{\X}(\frac{2\pi k}{Mq}+s_{pq}) + \frac{2\pi}{Mq} \tilde{\T}([\frac{2\pi  (k+1)}{Mq}+s_{pq}]^-), & k=0,1,\ldots,Mq-2,\\
\tilde{\X}(2\pi ) = \tilde{\X}(\frac{2\pi  (k)}{Mq}+s_{pq}) +   (\frac{2\pi}{Mq} - s_{pq} ) \tilde{\T}([\frac{2\pi  (k+1)}{Mq}+s_{pq}]^-), & k=Mq-1,
\end{cases}
\end{align*}

\noindent with $\tilde{\X}(0)=(0,0,0)$, and $s_{pq}=  (2\theta_0 p /Mq) \bmod  (2\pi/Mq) $. Note that $\tilde{\X}(0)$ and $\tilde{\X}(2\pi)$ do not correspond to a corner, but we compute them to define the correct rotation $\mathbf{L}_1$, which allows the polygonal curve to be aligned in such a way that the third component of $\tilde{\X}$ is parallel to the $z$-axis. The matrix $\mathbf{L}_1$ performs a rotation of angle equal to the one between $\mathbf{v}= \tilde{\X}(2\pi)-\tilde{\X}(0)$ and $(0,0,1)$, about an axis orthogonal to a plane spanned by these two vectors. As a result, $\T = \mathbf{L}_1 \cdot \tilde{\T}$, $\X = \mathbf{L}_1 \cdot \tilde{\X}$, and we obtain $\X$ and $\T$, up to a rotation and a vertical movement. The vertical movement can be computed at any time $\tpq$ from the speed of the center of mass $c_M$; but the rotation appears to be more involved, and, in fact, it turns out to be a multifractal. We will discuss this further in the next section.

\section{Numerical experiments}

\label{sec:Num-simu}

\begin{figure}[!htbp]\centering
	\includegraphics[width=0.449\textwidth, clip=true, align=t]{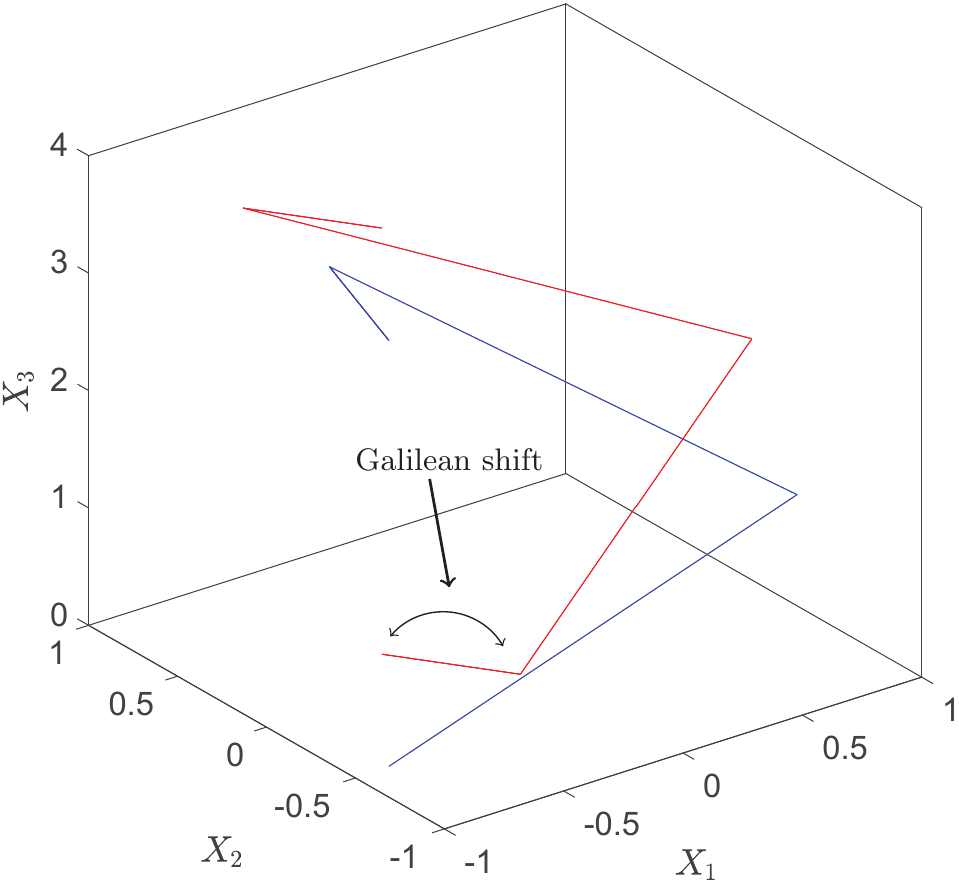}\includegraphics[width=0.551\textwidth, clip=true, align=t]{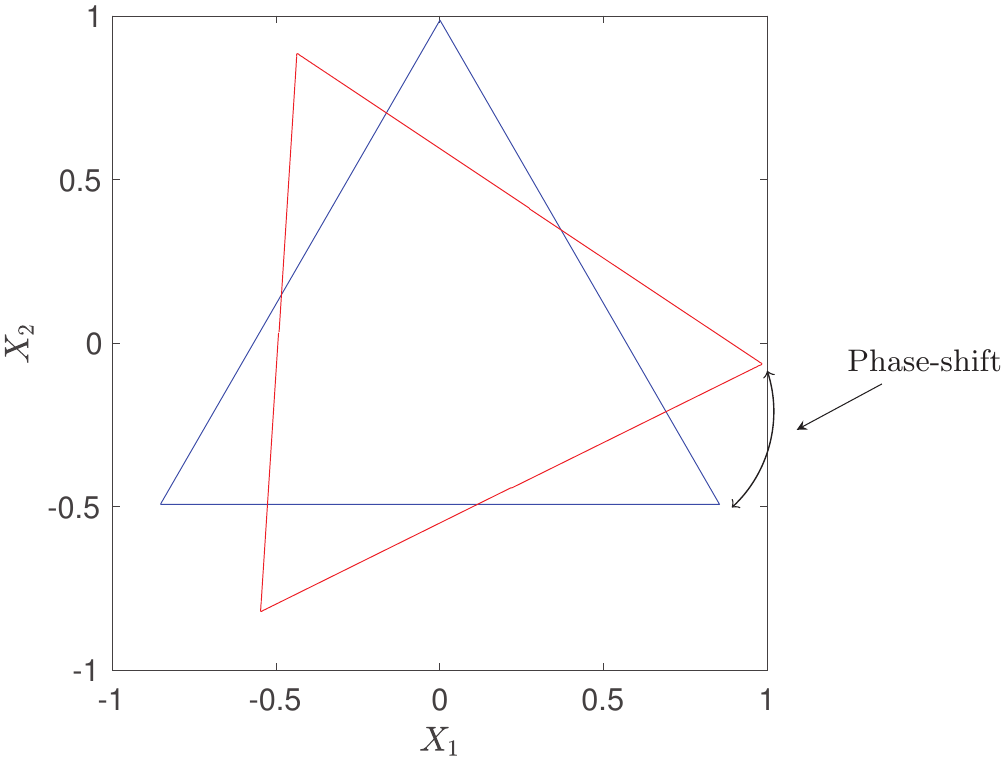}
	\caption{$\X(s,t)$, for $M=3$, $\theta_0=\pi/2,$ $b=0.5774\ldots$, at $t=0$ (blue), and at $t=T_f$ (red). The left-hand side shows the Galilean shift, and the right-hand side, the phase shift.} 
	\label{fig:M3-TWO-SHIFTS}
\end{figure}

We solve numerically VFE and the Schr\"{o}dinger map with the same numerical method as in \cite{HozVega2014} (see Section 4 and the references therein), to approximate \eqref{eq:VFE}-\eqref{eq:SMP} for the initial data $\X(s,0)$ and $\T(s,0)$, as given in (\ref{eq:Ini-data-T}) and (\ref{eq:Ini-data-X}), respectively. More precisely, the space interval $[0,2\pi)$ is discretized into $N$ equally spaced nodes, i.e. $s_k=2\pi k/N$, $k=0,1,\ldots,N-1$, and the time period $[0,T_f]$, with $T_f=2\pi/M^2$, into $N_t+1$ equally spaced time instants, taking $\Delta t = T_f / N_t$. Again, using the symmetries of $\T$, we can reduce all the discrete Fourier transforms of $N$ elements, to $N/M$ elements, reducing the computation cost quite effectively. With respect to the stability constraints on $N$ and $\Delta t$, they are also the same as in \cite{HozVega2014}; hence, $N/M$ and $N_t$, once fixed for one value of $M$, can be used for all $M$. On the other hand, since $\Delta t = \mathcal{O}(1/M^2)$, we can expect more accurate results for larger $M$. In our numerical simulations, we have taken $N/M= 480\cdot 2^r$, $N_t = 136080\cdot 4^r$, $r=0,1,\ldots$, and different values of $b$ (or $\theta_0$).

Recall that the initial curve is characterized by two parameters, $M$ and $b$. When $b\in(0,1)$, as $M$ is increased, the resulting initial curve tends to a helix. On the other hand, for a fixed $M$, as $b$ tends to $1$, the curve approaches a straight line. In our numerical simulations, we have analyzed both limits and computed the relevant quantities in each case. Apart from the fact that, at any rational time $\tpq$, there are $Mq$ or $Mq/2$ corners in $s\in[0,2\pi)$ (depending on whether $q$ is odd or even), we observe that the evolution is not periodic in time. As mentioned above, due to the Galilean shift, a corner initially located at $s=0$ moves to $s=2\theta_0/M$ at the end of one time-period; and the new helical $M$-polygon is rotated counterclockwise with respect to the $z$-axis by a certain amount, which we refer to as the phase shift. Figure \ref{fig:M3-TWO-SHIFTS} shows both shifts, for $M=3$ and $\theta=\pi/2$.

\begin{figure}[!htbp]\centering
	\includegraphics[width=0.5\textwidth, clip=true]{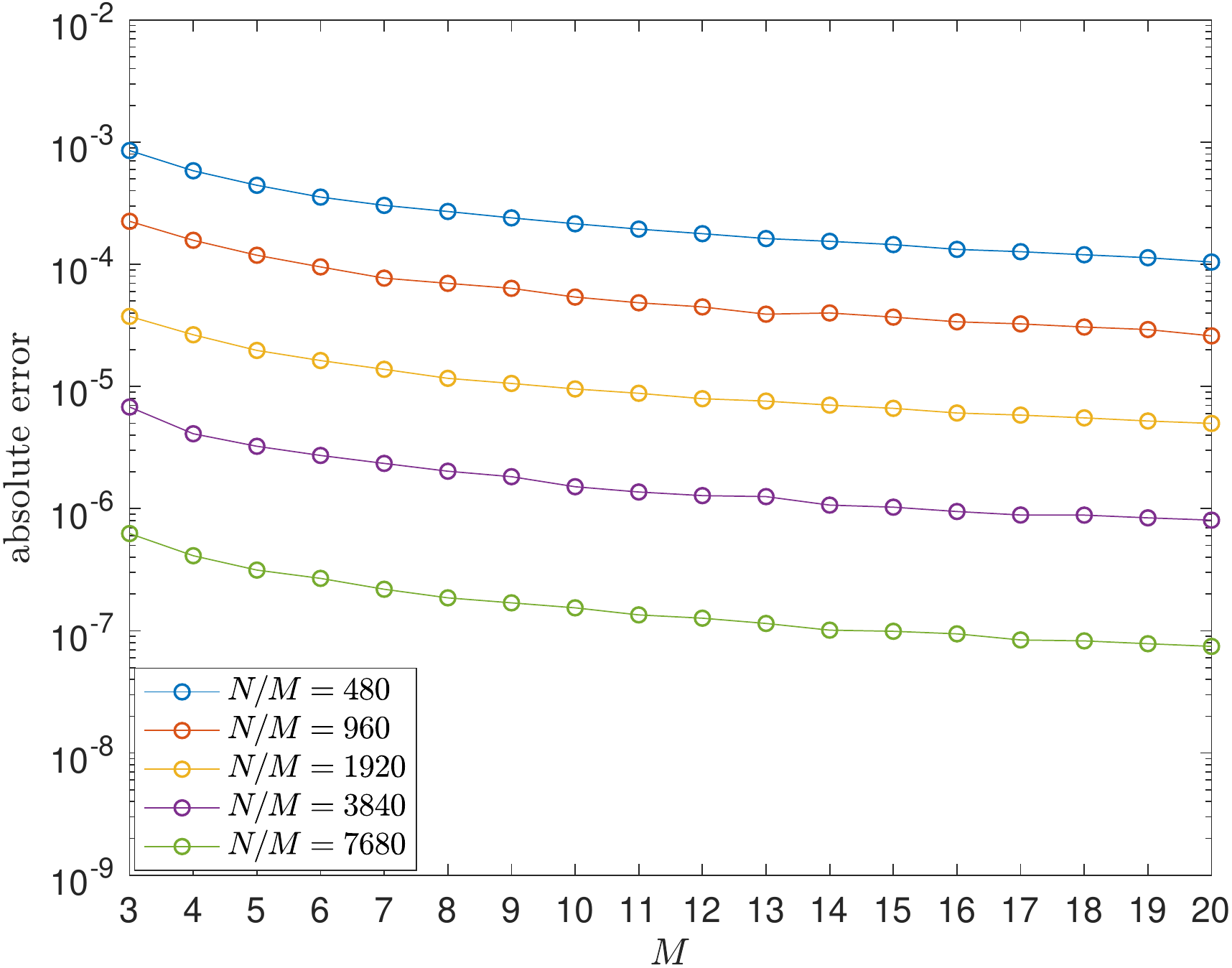}\includegraphics[width=0.5\textwidth, clip=true]{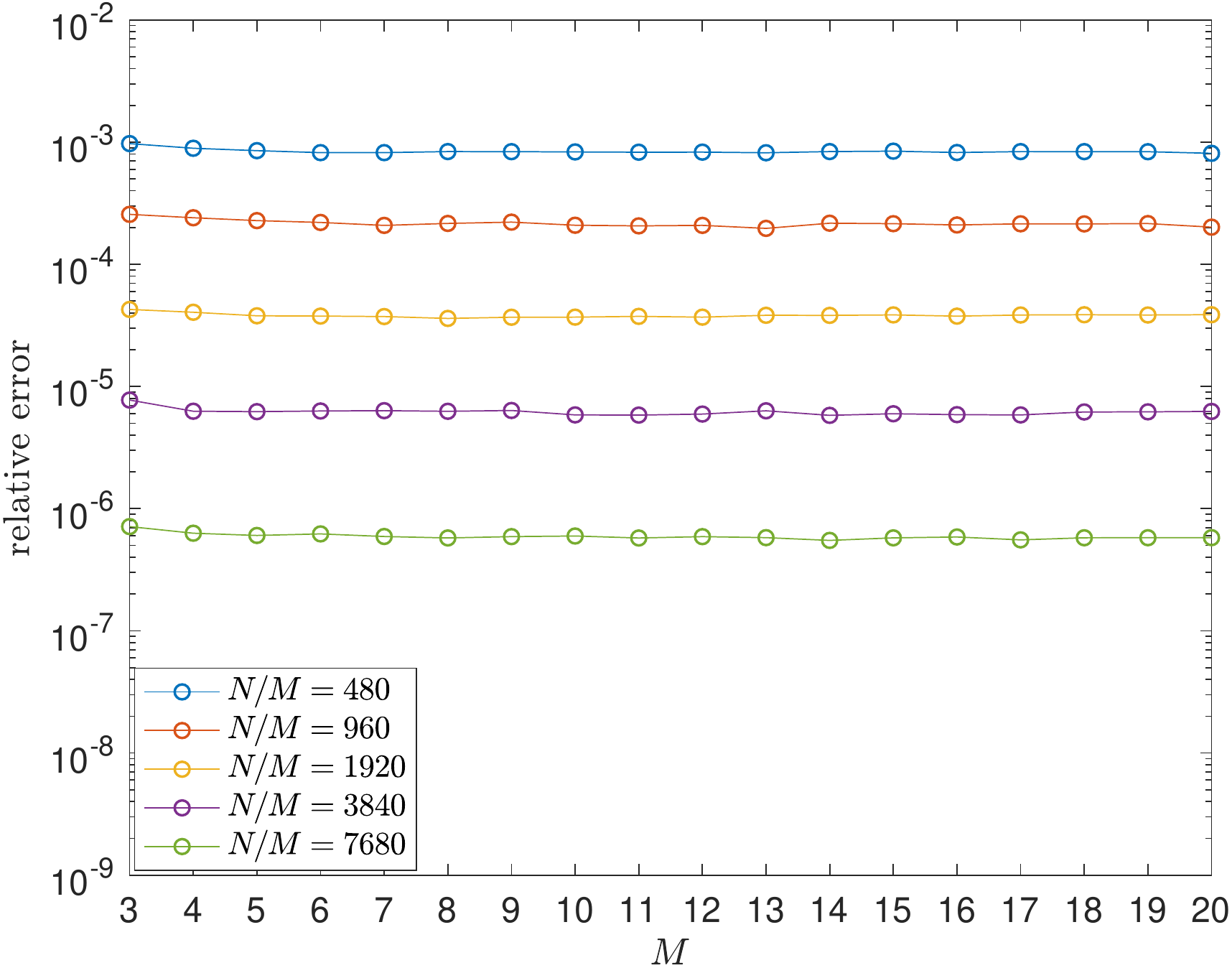}
	\caption{Absolute error (left) and relative error (right) as in (\ref{eq:err-rhoq}), in semilogarithmic scale, for the angle $\rho_q$, taking $b=0.4$, $q=5$, $M=3,4,\ldots,20$, and different values of $N/M$. The error clearly decreases, as $N/M$ increases, showing the convergence between the numerical and theoretical values.}
	\label{fig:rho-q5-b04}
\end{figure}

It is also important to compare the numerical value of $\rho_q$ with the one given by the algebraic expression in (\ref{eq:rho_q}), so for a given $q$, we compute the following errors:
\begin{equation}\label{eq:err-rhoq}
\left\{
	\begin{aligned}
	\Delta \rho_{q,N/M,M}^{abs} &=\max\limits_{\stackrel{p\in\{0,1,\ldots,q-1\}}{\gcd(p,q)=1}}\max\limits_{j=0,1,\ldots,Mq-1} \left|\rho_q-\rho_{pq}^{num,j}\right|,  \\
	\Delta \rho_{q,N/M,M}^{rel}&=\max\limits_{\stackrel{p\in\{0,1,\ldots,q-1\}}{\gcd(p,q)=1}}\max\limits_{j=0,1,\ldots,Mq-1} \left|\frac{\rho_q-\rho_{pq}^{num,j}}{\rho_q}\right|,	
	\end{aligned}
\right.
\end{equation}

\noindent where $\rho_{pq}^{num,j} = \arccos (\T_j \cdot \T_{j+1})$, $j=0,1,\ldots,Mq-1$. The value of the tangent vectors $\T_j$, which are piecewise constant at every time $\tpq$, has been calculated using the mean of the inner points; for example, for $\T(s)$, with $s\in[0,2\pi/Mq)$, we take the mean of the values corresponding to $s\in [\pi/2Mq, 3\pi/2Mq)$, etc. Using (\ref{eq:err-rhoq}), we compute the absolute and relative errors, for different values of $b$, $q$ and $M$. The results for $b=0.4$, $q=5$, $M=3,4,\ldots,20$, $N/M=480,960,\ldots,7680$, are plotted in Figure \ref{fig:rho-q5-b04}. Note that each color corresponds to a different discretization in the numerical scheme, which clearly shows the convergence of the errors, and hence, the agreement between the numerical and the algebraic values.

Besides the formation of new corners, it is evident from the numerical simulations that the evolution of a helical $M$-polygon involves a vertical upward movement as well. In order to determine the height of the polygonal curve at any time $t>0$, we compute its center of mass by taking the mean of all the values of $\X(s_j,t)$, i.e.,
$h_{N,k}(t) = \frac{1}{N}\sum_{j=0}^{N-1} X_k(s_j,t)$, $k=1,2,3$. $h_{N,1}(t)$ and $h_{N,2}(t)$ are equal to zero, whereas $h_{N,3}(t)$ grows linearly with $t$. Thanks to the symmetries mentioned in (\ref{sec:Spatial_Sym_X_T}), we can obtain $h_{N,3}(t)$ by using only $N/M$ values, which gives $h_{N,3}(t) = h_{N/M,3}(t) + \pi b (M-1) / M$. Then, after removing the subscripts and denoting the height $h_{N,3}(t)$ as $h(t)$, the speed of the center of mass $c_M^{num}$ can be approximated numerically as $c_M^{num} = (h(2\pi/M^2) - h(0)) / (2\pi/M^2)$. On the other hand, by taking advantage of the spatial periodicity of $\T$, and using the approach described in \cite[Section 4]{HozVega2018}, the exact value of $c_M$ can be deduced algebraically:
\begin{equation} \label{eq:c_M_alg}
c_M = \frac{-2\ln \cos (\rho_0/2)}{(\pi/M) \tan(\pi/M)} = \frac{\ln (1+\tan^2 (\rho_0/2))}{(\pi/M) \tan(\pi/M)}.
\end{equation}

\begin{figure}[!htbp]\centering
	\includegraphics[width=0.5\textwidth, clip=true]{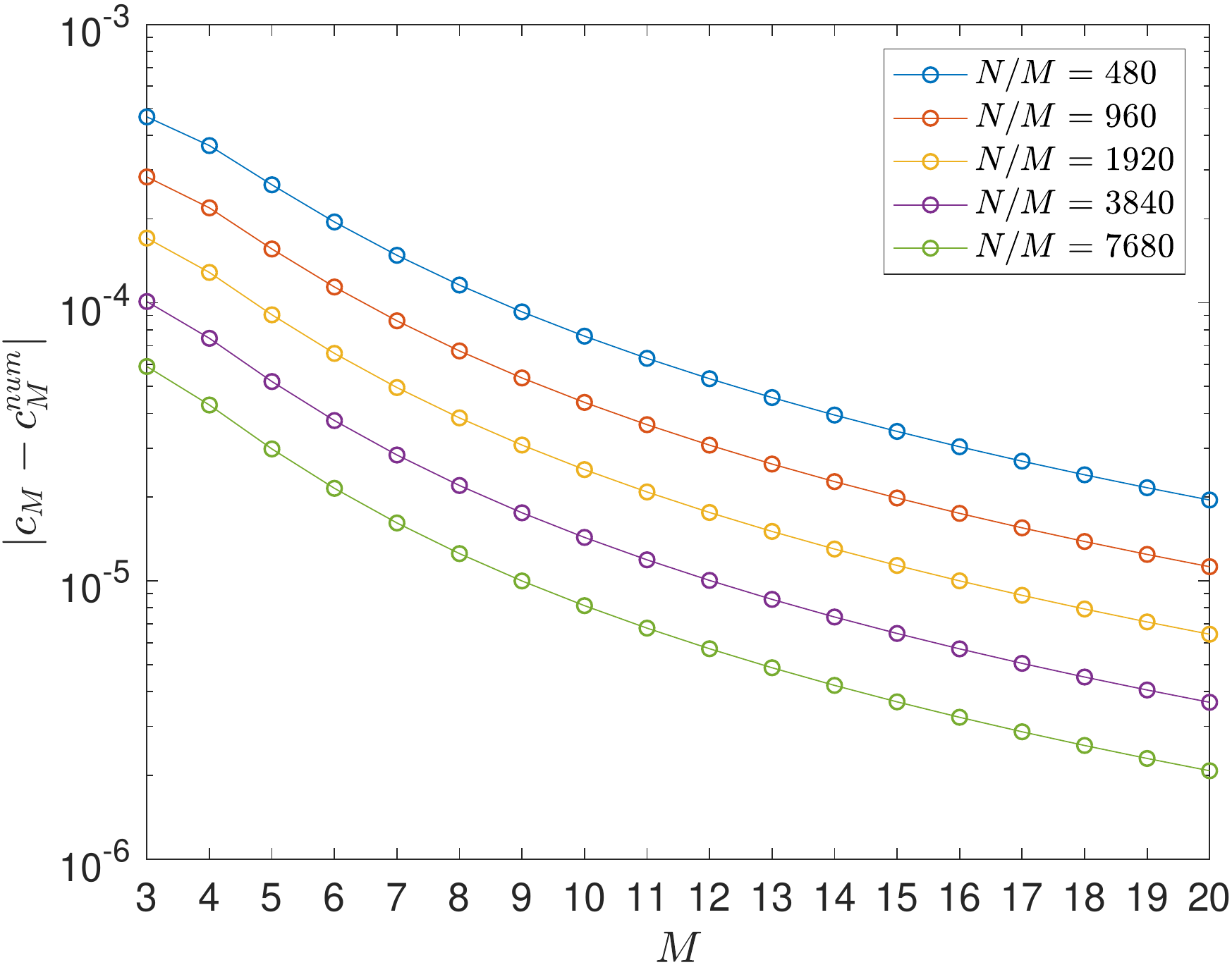}\includegraphics[width=0.5\textwidth, clip=true]{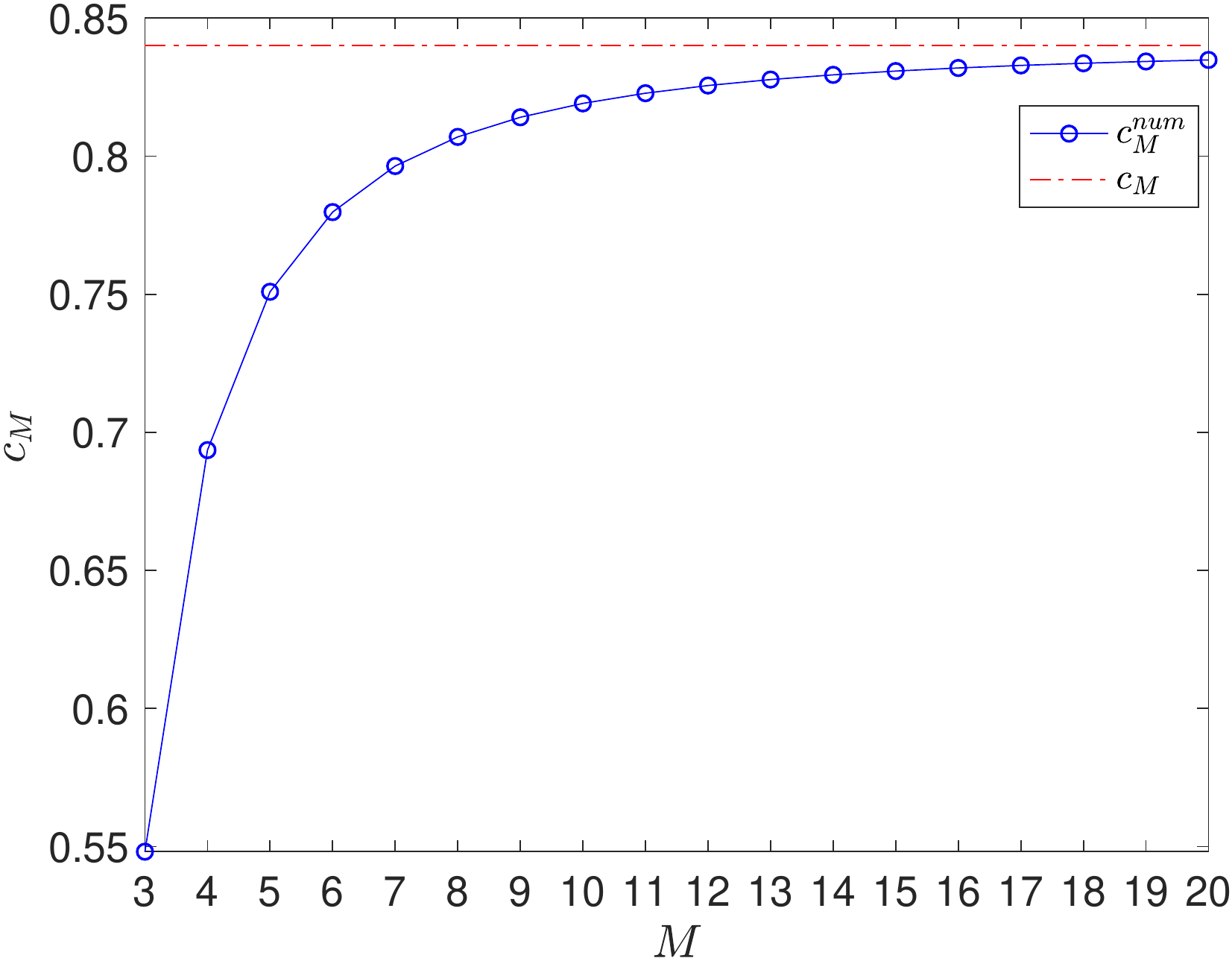}
	\caption{Left: $|c_M - c_M^{num}|$, computed for $b = 0.4$, and different values of $M$ and $N/M$. The error clearly decreases, as $N/M$ increases, showing the convergence between the numerical and the theoretical values. Right: $c_M$, computed using (\ref{eq:c_M_alg}), for different values of $M$. We have also plotted in dash-dotted line the limiting value $1-b^2=0.84$, to which we conjecture $c_M$ to tend, as $M\to\infty$.}
	\label{fig:cM-error}
\end{figure}

\noindent Remark that $\lim_{M\to \infty} c_M = 1-b^2$. Figure \ref{fig:cM-error} shows the value $|c_M - c_M^{num} |$, for $M=3,4,\ldots,20$, $b=0.4$, $N/M=480 \cdot 2^r$, $N_t = 136080 \cdot 4^r$, with $r=0,1,2,3,4$. We note that, for a given value of $M$, when $N/M$ is approximately doubled, the errors are divided by a factor slightly smaller than two, which implies that the errors behave as $\mathcal{O}((N/M)^{-1})$, and shows the convergence of $c_M^{num}$ to $c_M$. Moreover, as $M$ increases, the errors reduce, which can be explained from the fact that $\Delta t = \mathcal O(1/M^2)$. Finally, let us mention that $c_M^{num}$ converges to $1-b^2=0.84$, as $M\to\infty$.

In the following section, we will do a detailed study of $\X(0,t)$, which constitutes the main difference between the zero-torsion and the nonzero-torsion cases.

\section{Trajectory of $\X(0,t)$} \label{sec:X0t}

In the zero-torsion case, due to the symmetries of the closed $M$-polygons, the trajectory of one corner, i.e., $\X(0,t)$, which at the numerical level was claimed to be a multifractal, lies in a plane  \cite{HozVega2014}. However, for $\theta_0>0$, $\X(0,t)$ is no longer planar, and taking $t\in[0,T_f]$ is not enough to understand its structure, so we consider larger times multiple of $T_f$. Figure \ref{fig:X0t-M6} corresponds to an $M$-polygon with $M=6$, $\theta_0=\pi/5$, i.e., $b=0.5628\ldots$, $t\in[0,10\pi/3]$. Observe that $\X(0,t)$ has a helical shape, and exhibits a conspicuous fractal structure that repeats periodically, with some rotation and a vertical movement. In order to further understand this curve, we analyze each component of $\X(0,t)$, by using a Fourier series. In what follows, we consider two different cases, depending on the choice of the parameter $b$.

\subsection{Case with $b\in(0,1)$}\label{0<b<1}

\begin{figure}[!htbp]\centering
	\includegraphics[width=0.5\textwidth, clip=true]{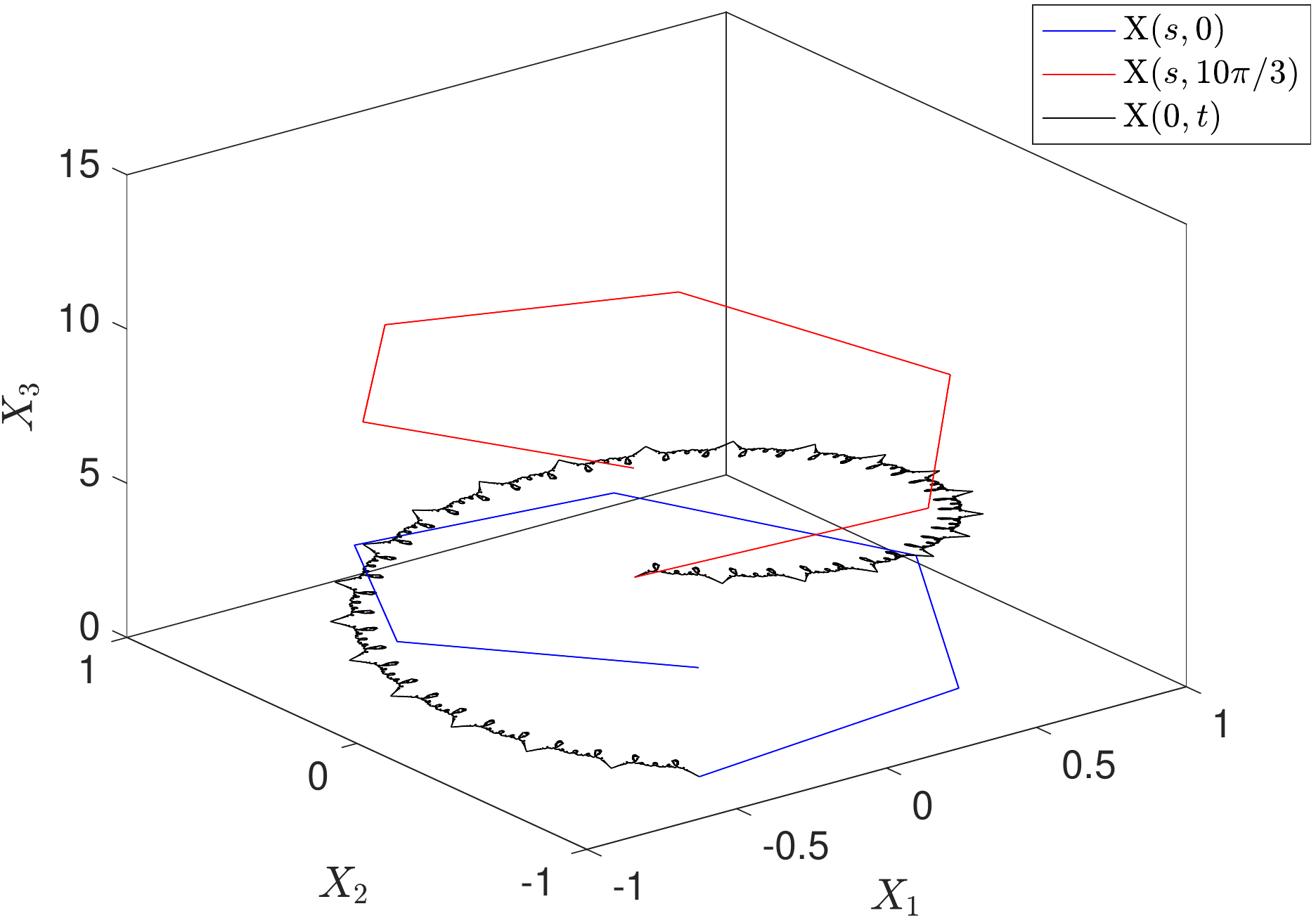}
	\caption{Initial polygon $\X(s,0)$ (blue), for $M=6$, $\theta_0=\pi/5$, i.e., $b=0.5628\ldots$, and its evolution at time $t =  10\pi/3$ (red), together with the curve described by $\X(0,t)$, for $t\in[0,10\pi/3]$ (black). $\X(0,t)$ has a
		conspicuous fractal structure which repeats periodically, with some rotation and a vertical movement.}\label{fig:X0t-M6}
\end{figure}

\begin{figure}[!htbp]\centering
	\includegraphics[width=0.327\textwidth, clip=true]{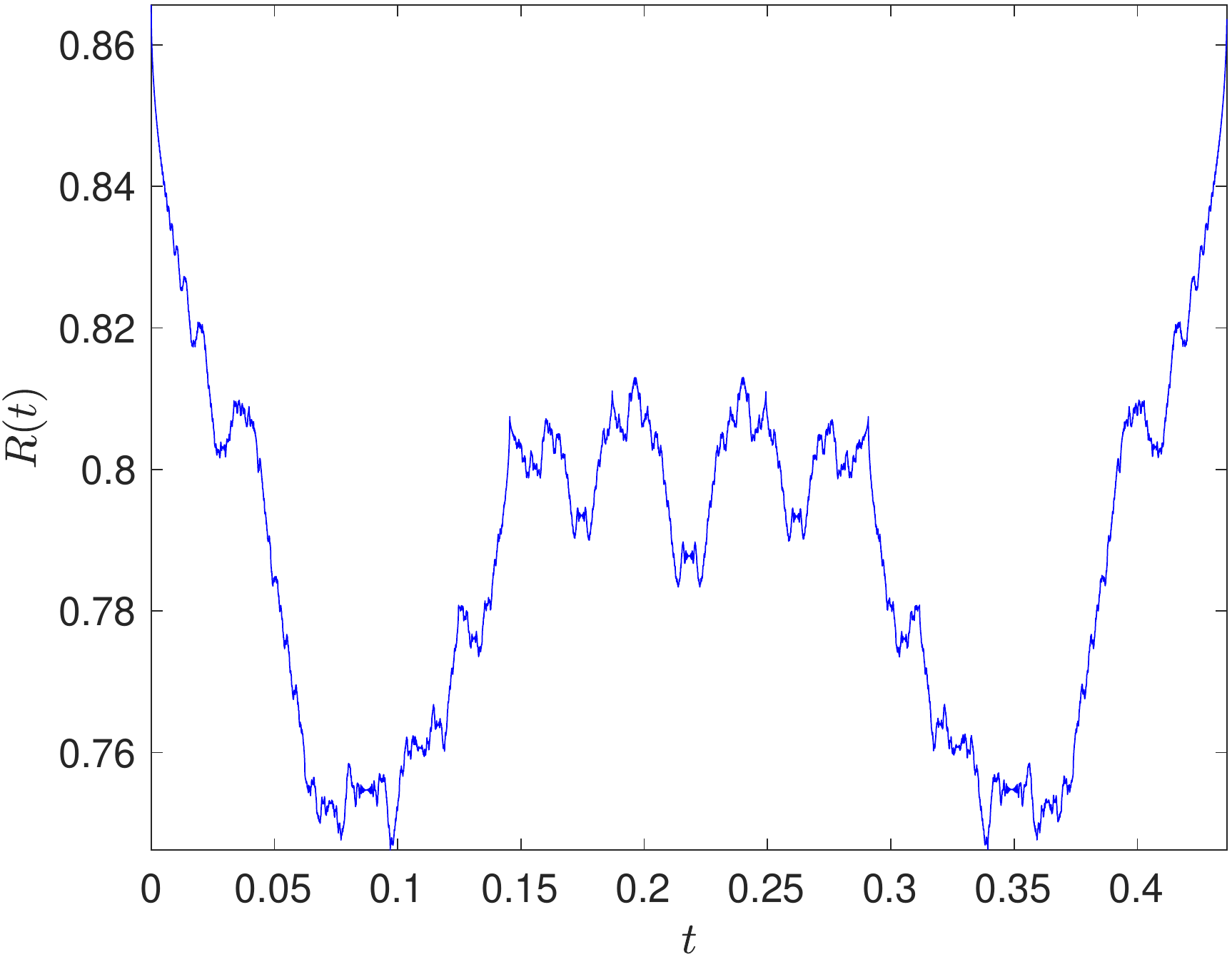}\includegraphics[width=0.346\textwidth, clip=true]{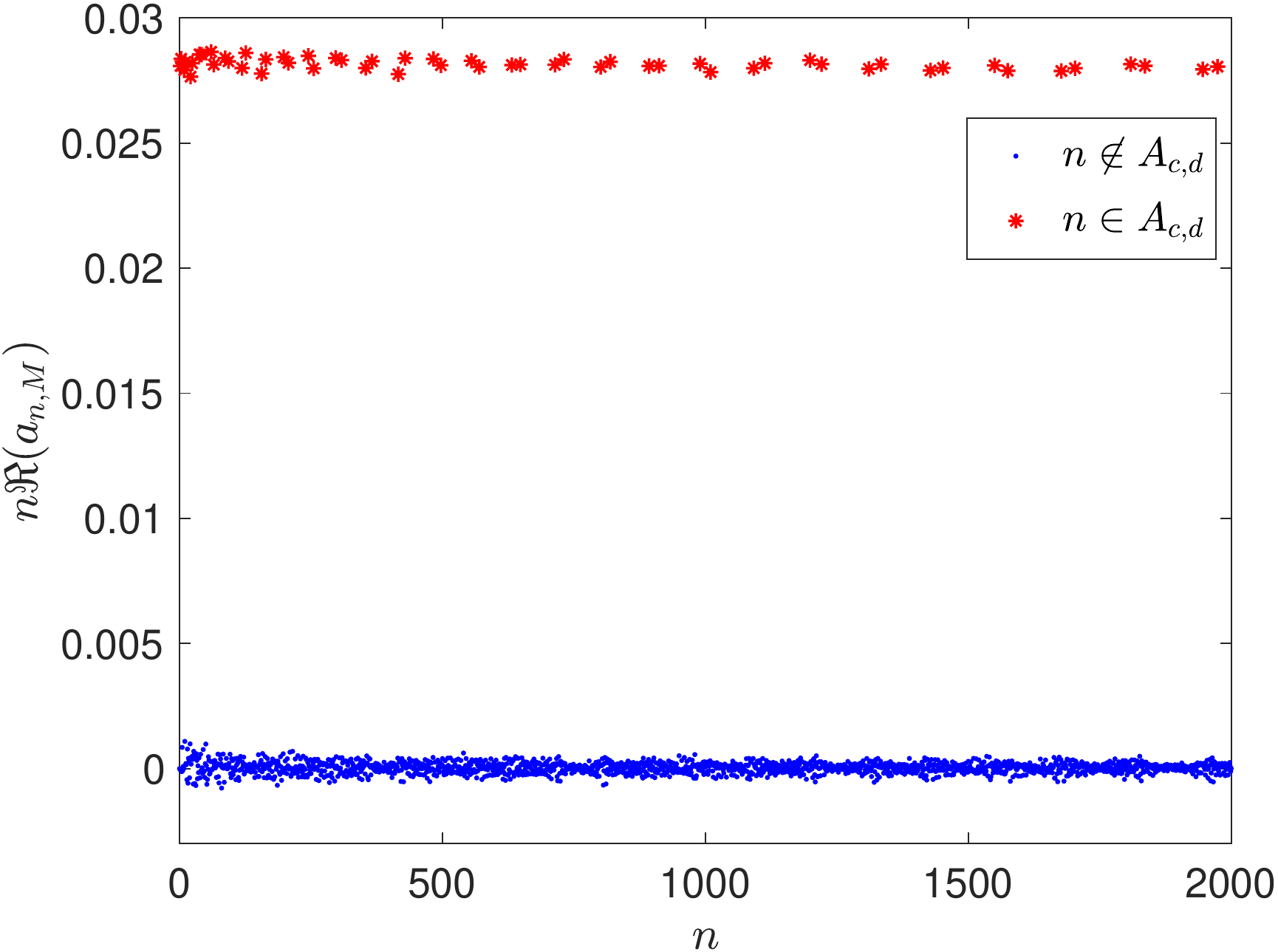}\includegraphics[width=0.327\textwidth, clip=true]{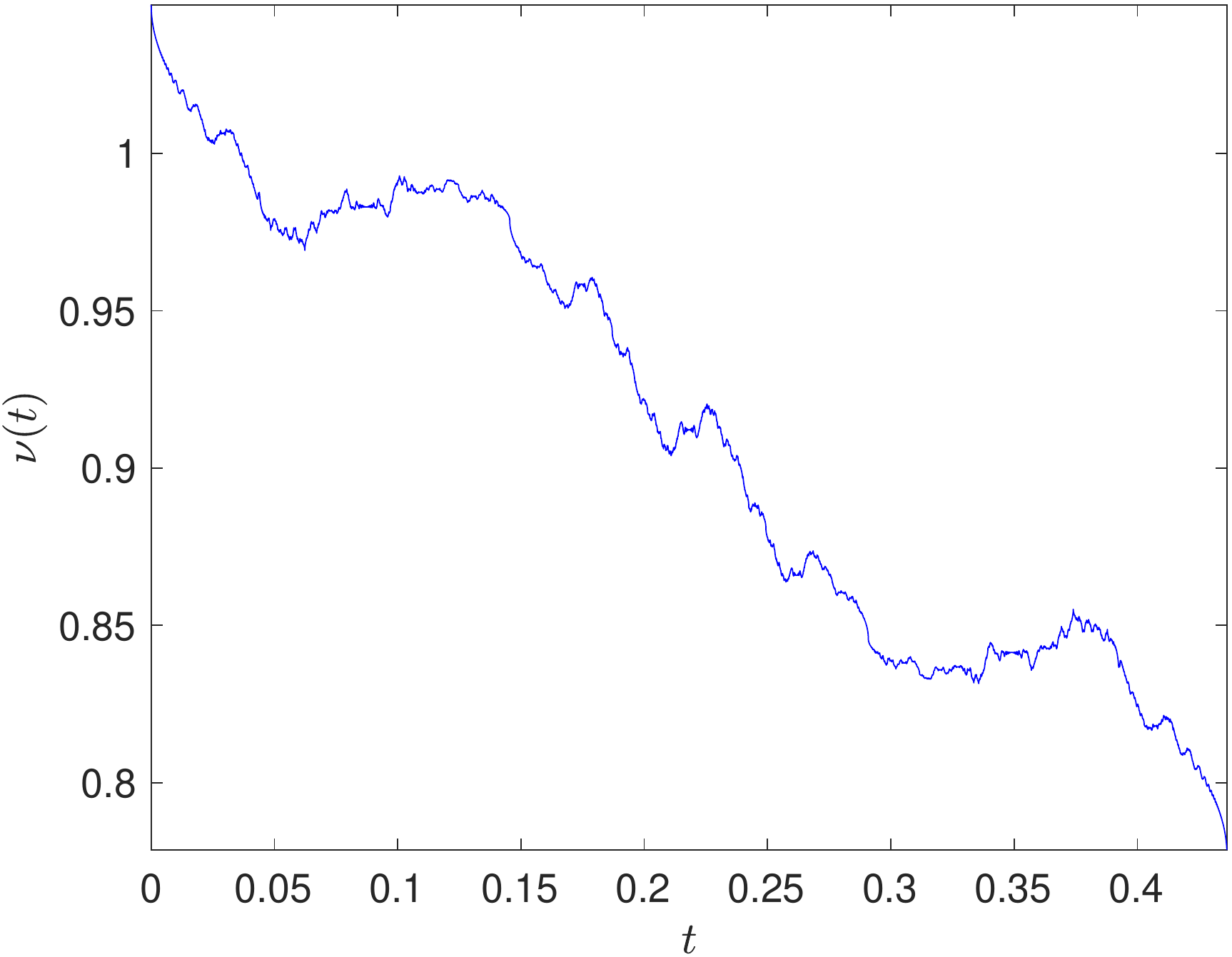}
	\caption{Left: $R(t)$, for $t\in[0,T_f^{c,d}]$, $M=6$, $\theta_0=\pi/5$, $c=1$, $d=5$, $b=0.5628\ldots$. Center: Approximation of $n\Re(a_{n,M})$, for $n = 1,2,\ldots,2000$, where $a_{n,M}$ are the Fourier coefficients of $R(t)$. The dominating points (red starred) are given by (\ref{eq:dom-pts}). Right: $\nu(t)$, which seems to have a multifractal structure as well.}\label{fig:Rt-M6}
\end{figure}

In the case of a planar $M$-polygon, the curve $\X(0,t)$ has a corner at those values of $t$, at which the polygonal curve $\X(s,t)$ has a corner at $s = 0$, i.e., at $t = t_{pq}$, with $q\not\equiv 2 \bmod 4$. However, when $\theta_0>0$, $\X(0,t)$ having corners depends mainly on the Galilean shift, or in other words, in order for the curve $\X(0,t)$ to have corners in finite time, $b$ should be chosen in such a way that, at $t = T_f$, the Galilean shift is a rational multiple of the side length $2\pi/M$ of the corresponding polygonal curve. This can be enforced by taking $\theta_0 = \pi c/d$, with $\gcd(c,d)=1$, $c,d \in \mathbb{N}$. Then, defining the following multiples of $T_f$:
\begin{equation}\label{eq:d_time_period}
\tcd \equiv
\begin{cases}
(d/2) T_f \equiv \pi d/M^2, & \mbox{if $c\cdot d$ odd,} \\
d\, T_f \equiv 2\pi d/M^2,  & \mbox{if $c\cdot d$ even,}
\end{cases}
\end{equation}

\noindent the numerical experiments reveal that, at times that are integer multiples of $\tcd$, $\X(0,t)$ has a corner, and the three-dimensional fractal structure of $X(0,t)$ repeats with period $\tcd$ (see Figure \ref{fig:X0t-M6})). On the other hand, at rational multiples of $\tcd$, $\X(0,t)$ has corners of different (smaller) scales. In order to better understand $\X(0,t)$, we define, in $t\in[0,\tcd]$:
$$
\begin{cases}
z_{1,2}(t) = X_1(0,t)+i X_2(0,t) = R(t) e^{i \nu(t)},
	\\
\tilde X_3(t) = X_3(0,t)-c_M t,
\end{cases}
$$

\noindent where $R(t) = \sqrt{X_1^2(0,t)+X_2^2(0,t)}$ and $\nu(t) =  \arctan( X_2(0,t)/X_1(0,t))$ give the polar representation of $z_{1,2}(t)$, and $\tilde X_3(t)$ is $X_3(0,t)$ without its vertical height. Since $R(t)$ and $\tilde X_3(t)$ are periodic, we can consider their Fourier expansion:
\begin{align*}
R(t) &= \sum\limits_{n=-\infty}^{\infty} a_{n,M} e^{2\pi i\,n\,t / \tcd}, \quad  
\tilde X_3(t) = \sum\limits_{n=-\infty}^{\infty} b_{n,M} e^{2\pi i\,n\,t / \tcd}, \quad t \in [0,T_f^{c,d}].
\end{align*}
We have approximated the Fourier coefficients $a_{n,M}$ and $b_{n,M}$ using the \textsc{MATLAB} command \texttt{fft}, for $M=6$, $\theta_0=\pi/5$, $t\in[0,5\pi/18]$. In the left-hand side of Figure \ref{fig:Rt-M6}, we have plotted $R(t)$; and, in the center, the real part of the approximations of $n\,a_{n,M}$, for $n = 1,2, \ldots, 2000$. This last plot can be referred to as the fingerprint of $R(t)$, and it is a tool to understand the multifractal structure of $R(t)$ (see also \cite[Figure 11]{HozVega2018}). Finally, on the right-hand side, we have plotted $\nu(t)$, which describes the angular movement of $\X(0,t)$ in the XY-plane and can be associated with the phase shift corresponding to the angular movement of a corner initially located at $s=0$. From its definition, one can suspect that $\nu(t)$ has a multifractal structure, too; hence, computing the phase shift at any rational time appears to be involved and deserves further research. We make some more comments on this in the next section.

We have carried out a careful study of the fingerprints of $R(t)$ and $\tilde X_3(t)$, with $t \in [0,T_f^{c,d}]$, for many different values of $M$, $c$ and $d$, and have found strong evidence that, when $\theta_0=c \pi / d$, $\gcd(c,d)=1$, the dominating points of the fingerprints belong to the following set:
\begin{equation}\label{eq:dom-pts}
A_{c,d}=\begin{cases}
\{n(n\,d + c)/2\,|\, n \in\mathbb Z \} \cap \mathbb N, & \mbox{if $c\cdot d$ odd,}
\\
\{n(n\,d + c)\,|\, n \in\mathbb Z \} \cap \mathbb N, &  \mbox{if $c\cdot d$ even.} 
\end{cases}
\end{equation}
On the other hand, the relation between $\X(0,t)$ for planar $M$-polygons and Riemann's non-differentiable function suggests comparing $\tilde X_3(t)$ in the helical $M$-polygon case and the imaginary part of 
\begin{equation}\label{eq:phi-c-d}
\phi_{c,d}(t) =  \sum\limits_{k \in A_{c,d}} \frac{e^{2\pi i k t}}{k}, \quad
t\in
\begin{cases}
[0,1/2], & \mbox{if $c\cdot d$ odd,} \\
[0,1], & \mbox{if $c\cdot d$ even.} 
\end{cases}
\end{equation}

\noindent where $A_{c,d}$ is given by (\ref{eq:dom-pts}).
\begin{figure}[!htbp]\centering
	\includegraphics[width=0.33\textwidth, clip=true]{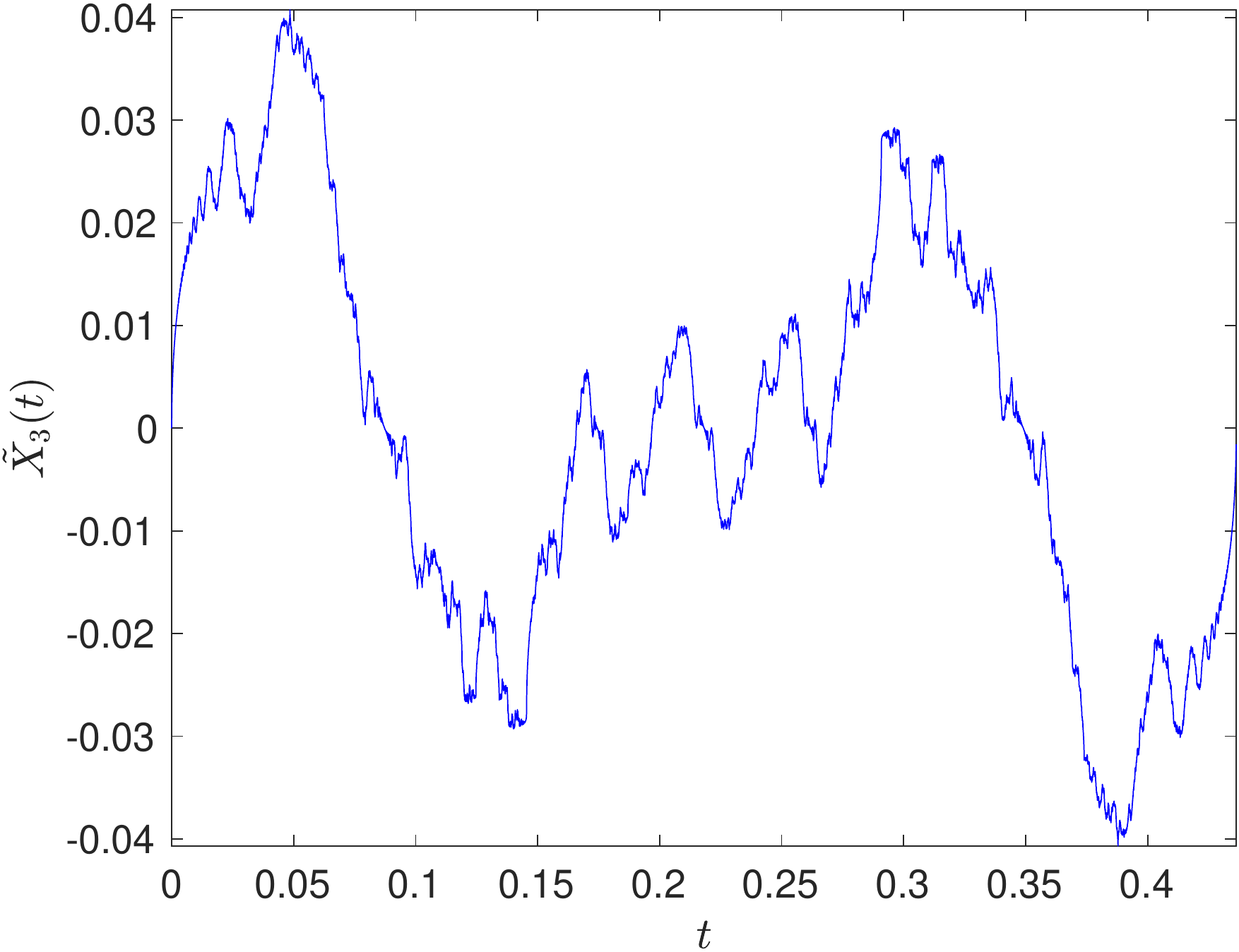}\includegraphics[width=0.34\textwidth, clip=true]{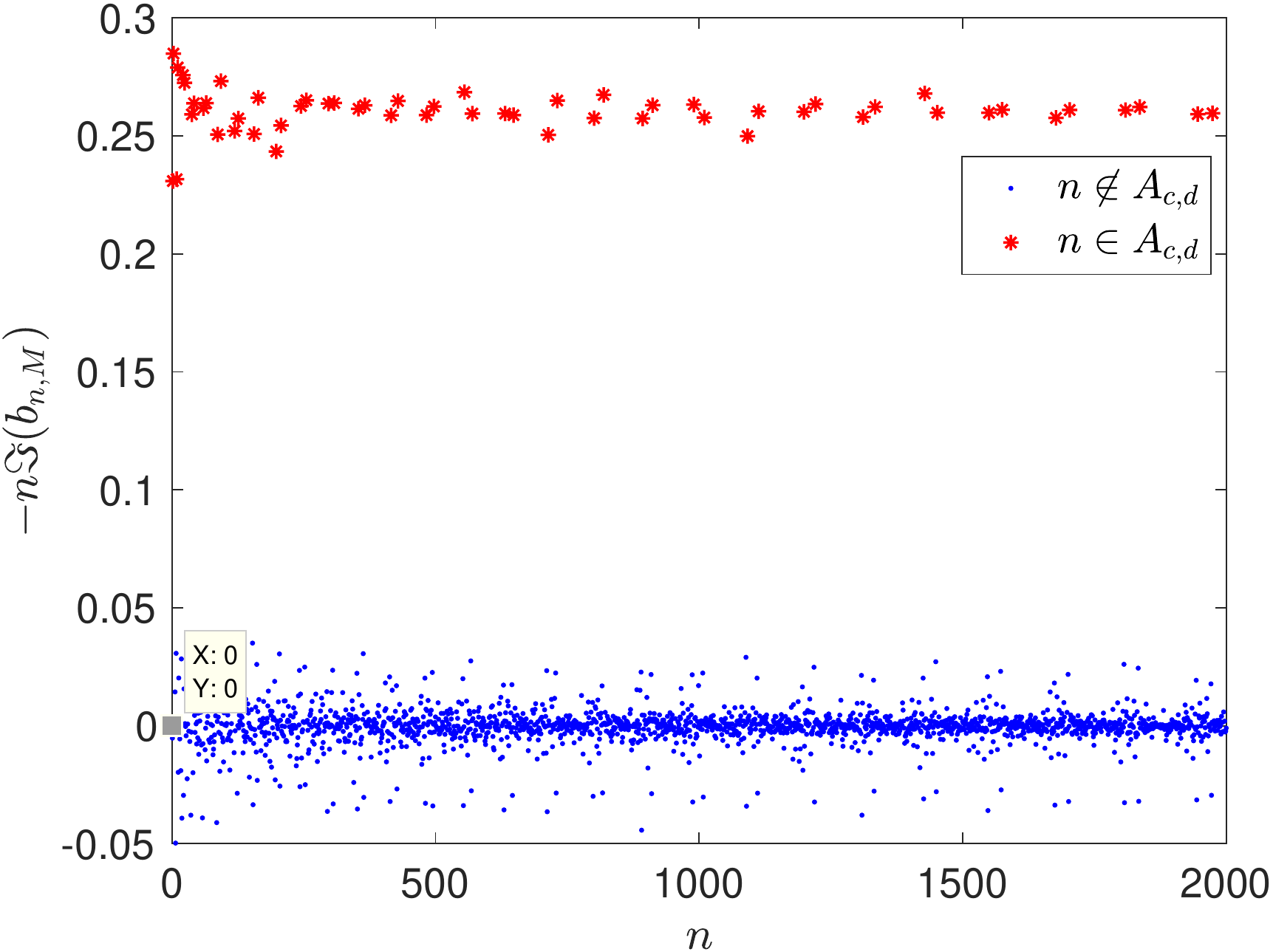}\includegraphics[width=0.33\textwidth, clip=true]{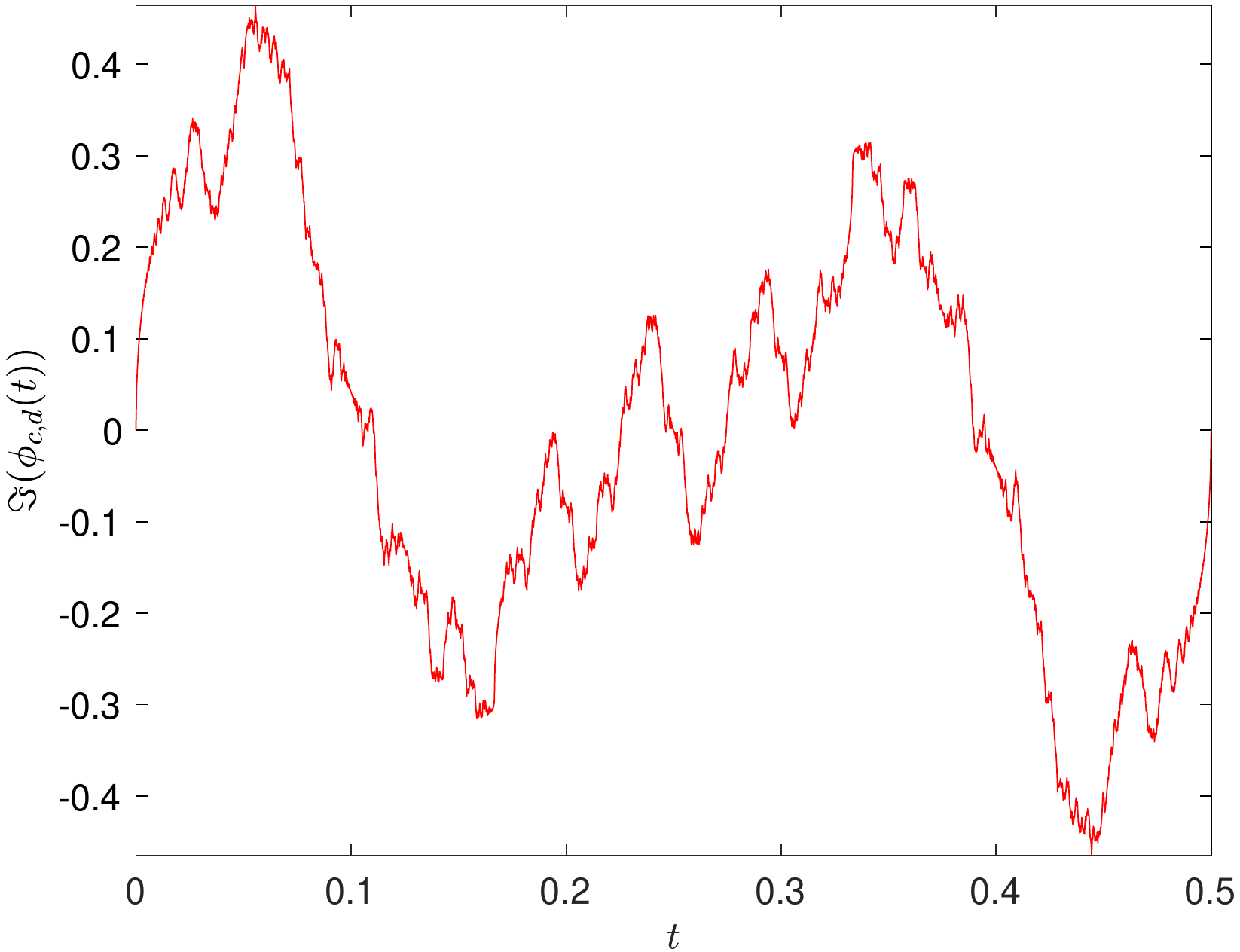}
	\caption{Left: $\tilde X_3(t)$, for $t\in[0,T_f^{c,d}]$, $M=6$, $\theta_0=\pi/5$, $c=1$, $d=5$, $b=0.5628\ldots$. Center: Approximation of $-n\Im(b_{n,M})$, for $n=1,2,\ldots, 2000$, where $b_{n,M}$ are the Fourier coefficients of $\tilde X_3(t)$ multiplied by $-1$. The dominating points (red starred) are given by (\ref{eq:dom-pts}). Right: Imaginary part of $\phi_{c,d}(t)$ in (\ref{eq:phi-c-d}), where the sum is taken over $2^{11}$ values.}\label{fig:XM3-M6}
\end{figure}

\noindent Figure \ref{fig:XM3-M6} is the continuation of Figure \ref{fig:Rt-M6}, and hence, all the parameters are identical. On the left-hand side, we have plotted $\tilde X_3(t)$; in the center, minus the imaginary part of the approximations of $n\,b_{n,M}$, for $n = 1,2, \ldots, 2000$, i.e., the fingerprint of $\tilde X_3(t)$; in general, the dominating points appear to be distributed around $1/4$ (when $c \cdot d$ is odd) or $1/2$ (when $c \cdot d$ is even). Finally, on the right-hand side, we have plotted the imaginary part of $\phi_{c,d}(t)$ which, except for a scaling, is visually indistinguishable from the $\tcd$-periodic curve $\tilde X_3(t)$ on the left-hand side.

\begin{figure}[!htbp]\centering
	\includegraphics[width=0.7\textwidth, clip=true]{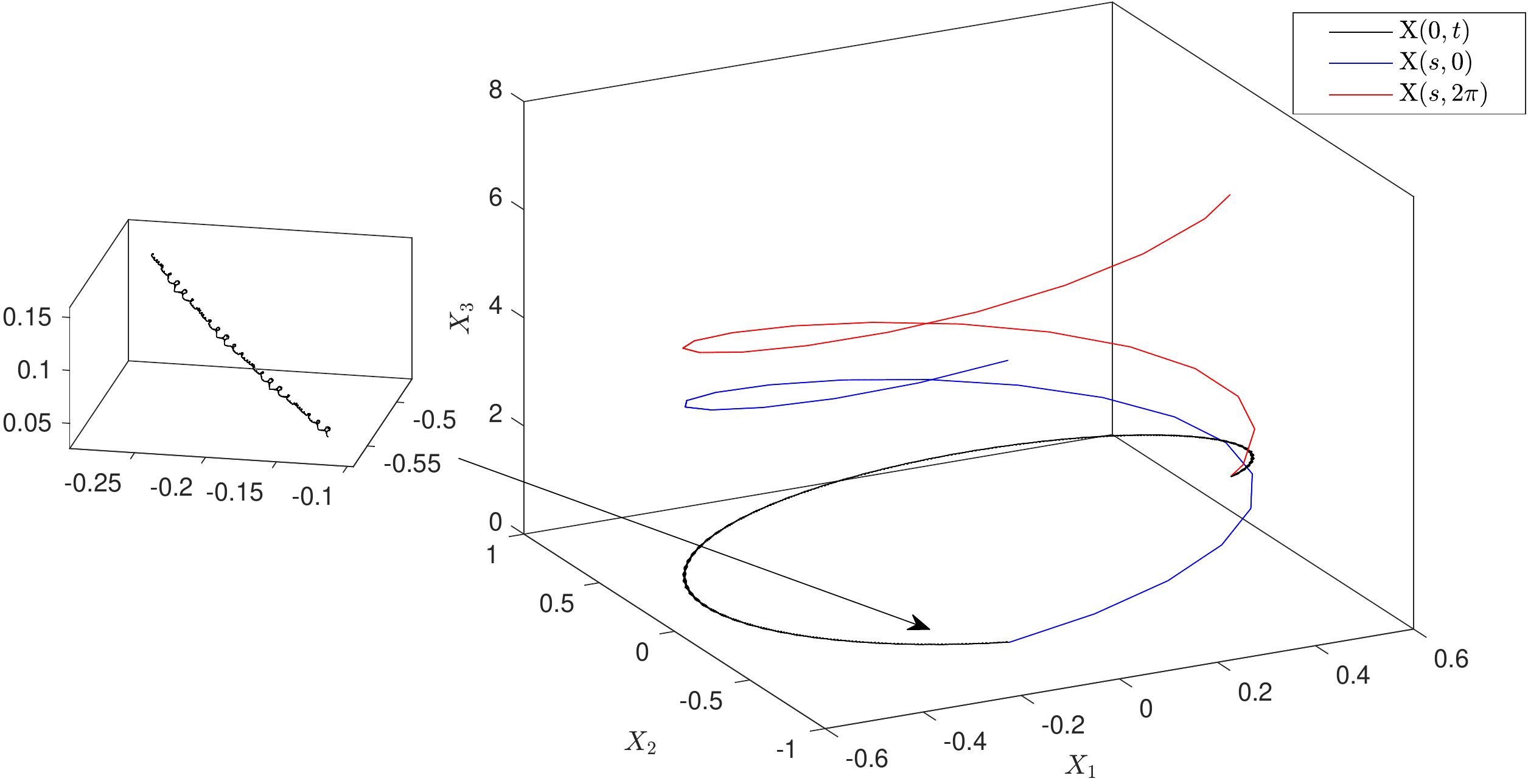}
	\caption{Initial polygon $\X(s,0)$ (blue), for $M=20$, $\theta_0=\pi/12$, i.e., $b=0.8312\ldots$, and its evolution at time $t  = 2\pi$ (red), together with the curve described by $\X(0,t)$, for $t\in[0,2\pi]$ (black). $\X(0,t)$ has a conspicuous fractal helical structure, as shown in the zoomed part.}\label{fig:M20-X0t}
\end{figure}

\subsubsection{Computation of $\lim_{M\to \infty} \X(0,t)$}

\label{sec:sub2-Rot-shift} 
\begin{figure}[!htbp]\centering
	\includegraphics[width=0.495\textwidth, clip=true]{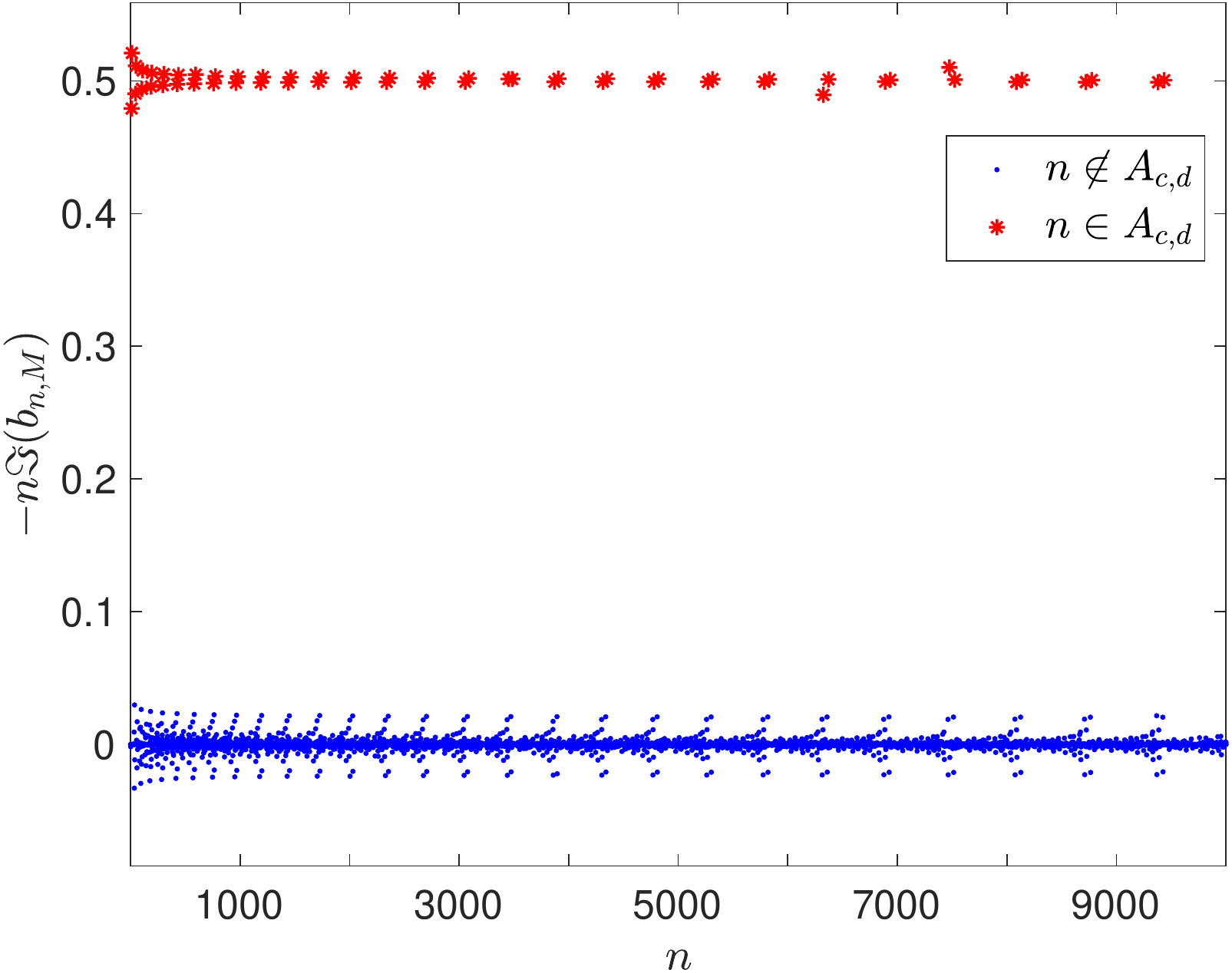}\includegraphics[width=0.505\textwidth, clip=true]{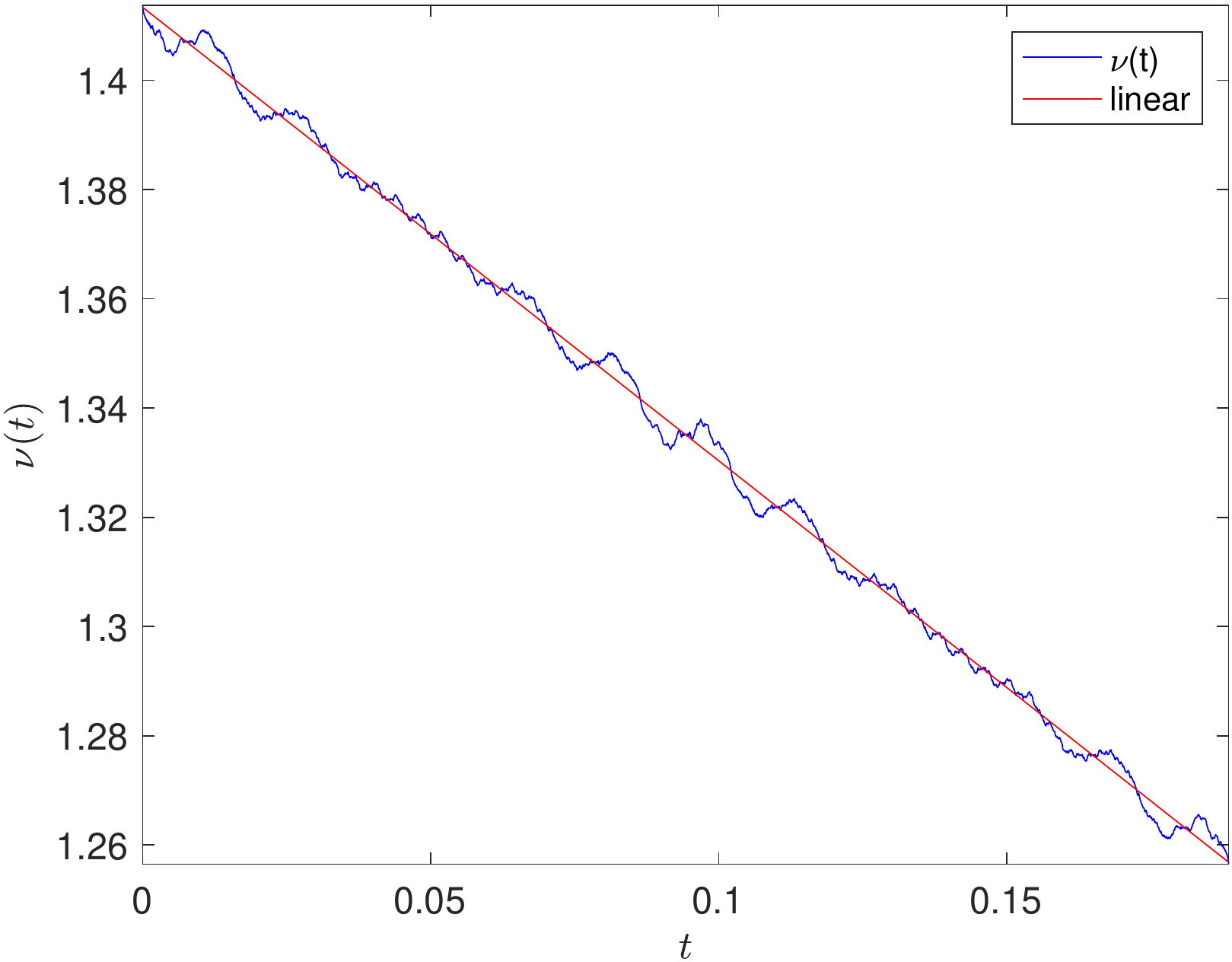}
	\caption{Left: fingerprint of $\tilde X_{3}(t)$, taking $M=20$, $\theta_0=\pi/12$, $c=1$, $d=12$, $b=0.8312\ldots$, $t\in [0,T_f^{c,d}]$. The convergence of the dominating points to $0.5$ is clearly visible. Right: $\nu(t)$, together with its best fitting line; the slope is $-0.8304 \ldots$, and its modulus is close to the value of $b$.}\label{fig:M20-XM3-nu}
\end{figure}

As we have seen, when $b\in(0,1)$, the curve $\X(0,t)$ is not periodic, but studying its structure componentwise sheds light on its behavior. We have also considered its time evolution for $M\gg1$ and sufficiently large values of $\tcd$. Figure \ref{fig:M20-X0t} shows $\X(0,t)$, for $M=20$, $\theta_0=\pi/12$, $b=0.8312\ldots, t\in[0,12\pi/5]$, along with the initial and final helical polygonal curve. The fingerprint of $\tilde X_{3}(t)$ has been plotted on the left-hand side of Figure \ref{fig:M20-XM3-nu}, and it is quite clear that the dominating points converge to $1/2$; indeed, the convergence is stronger than that in the center of Figure \ref{fig:XM3-M6}, since we have take $10^4$ points on the $x$-axis. This and other numerical experiments enable us to conjecture that 
\begin{equation*}
 \lim\limits_{M\to \infty} |n\, b_{n,M}|  = 
 \begin{cases}
1/4, & \text{if $n \in A_{c,d}$ and $c\cdot d$ odd},
\\
1/2, & \text{if $n \in A_{c,d}$ and $c\cdot d$ even},
\\
 0, & \text{otherwise}.
 \end{cases} 
 \end{equation*}

\begin{figure}[!htbp]\centering
	\includegraphics[width=0.503\textwidth, clip=true, align=t]{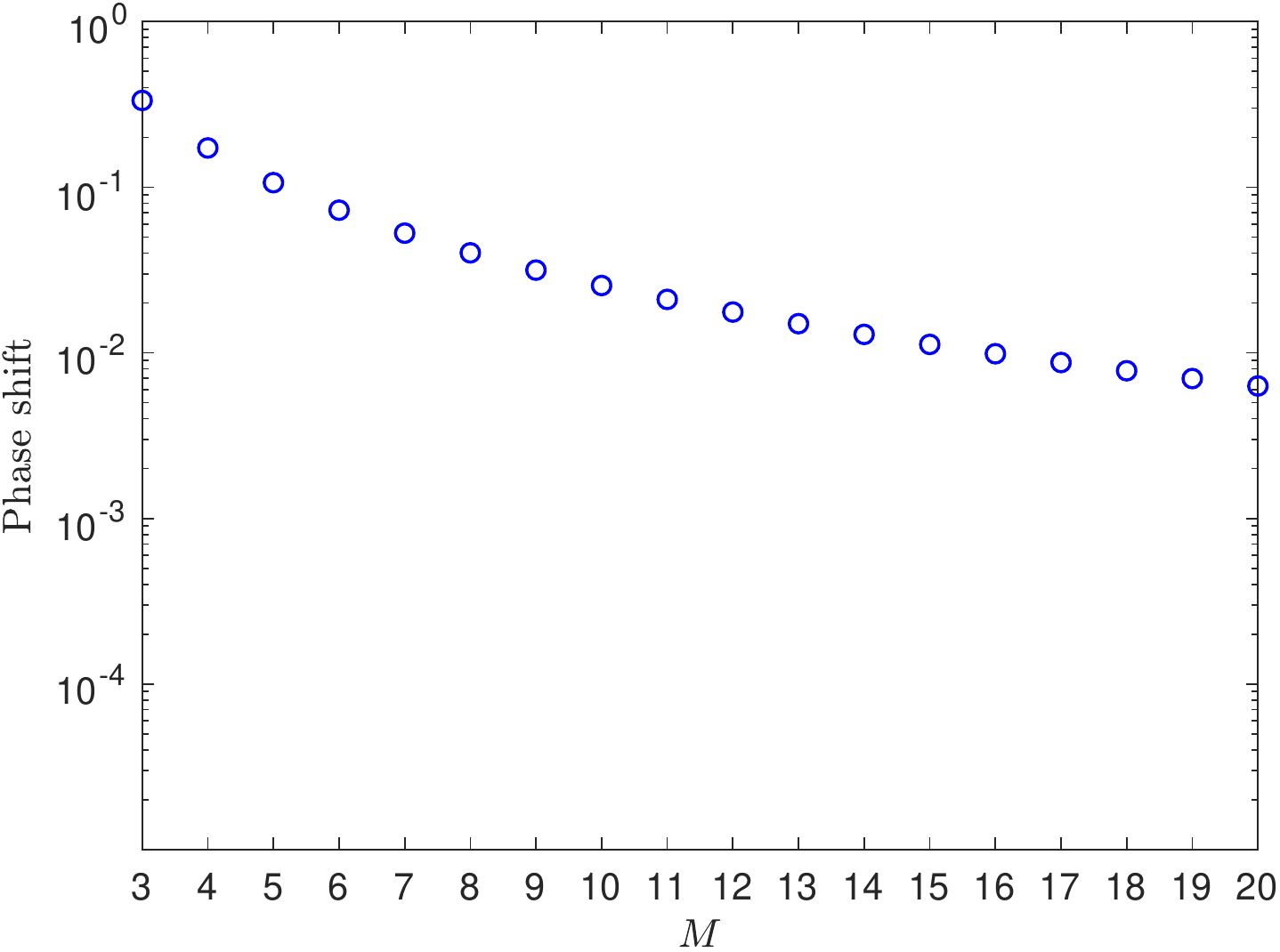}\includegraphics[width=0.497\textwidth, clip=true, align=t]{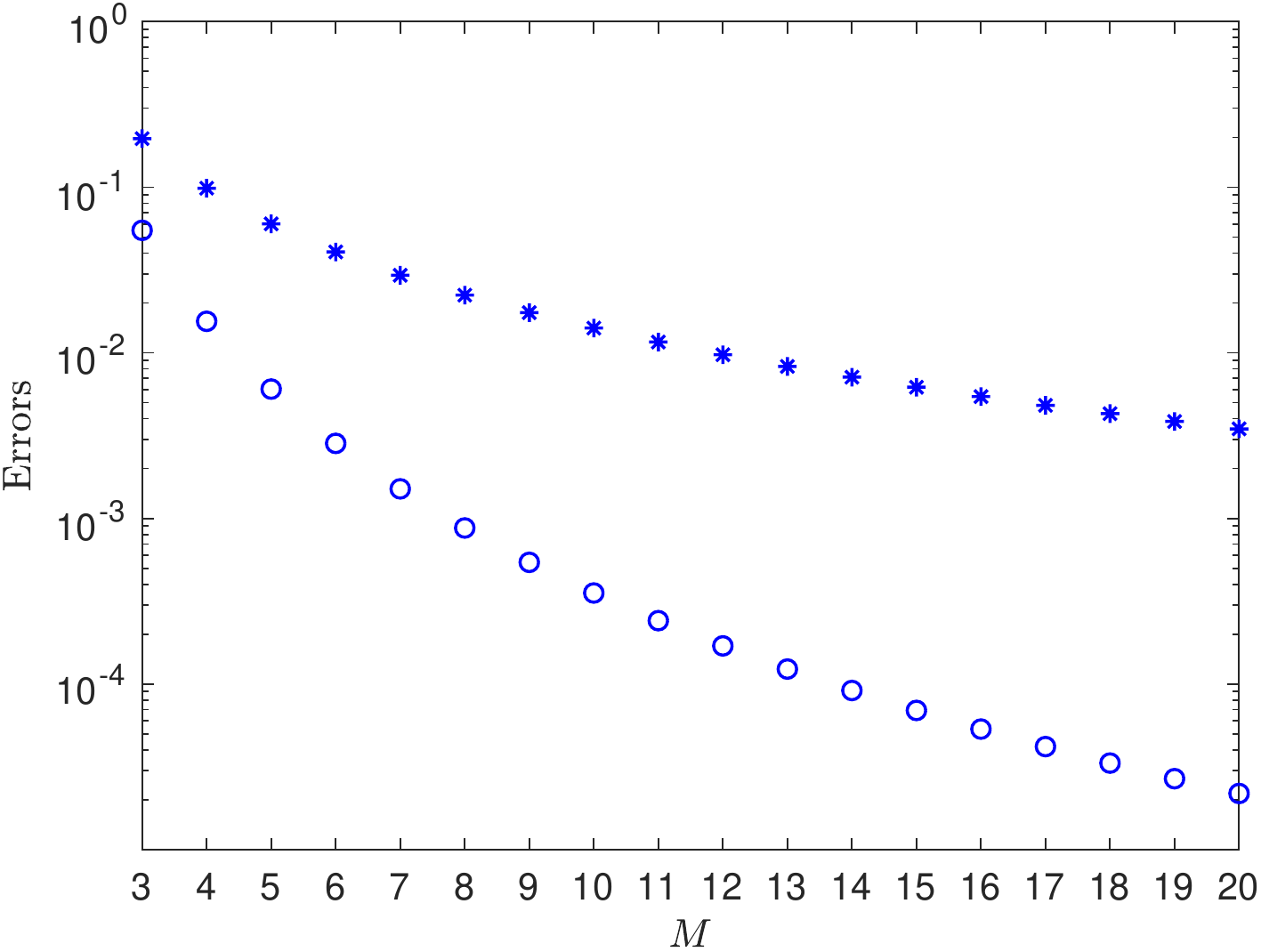}
	\caption{Left: Phase shift, for $M=3,4,\ldots,20$, $b=0.4$, $N/M=7680$. Right: absolute errors (circled points) and relative errors (starred points) computed by comparing the phase shift with $2\pi b/M^2$, for each $M$. For a given value of $b>0$, the phase shift at the end of one time-period decreases, as $M$ increases.}
	\label{fig:rot-shft-all}
\end{figure}

\noindent On the other hand, the curve $\nu(t)$ shown on the right-hand side of Figure \ref{fig:M20-XM3-nu} appears to converge to a straight line. In this regard, we have performed a basic linear fitting $\nu(t)=m\, t+c$ (red), with $c=1.4133\ldots$, $m=-0.8304\ldots$. Note that the modulus of $m$ can be compared with the value of $b$, which can also be regarded as the angular velocity of the helical curve $\X$, as it evolves in time. Continuing the discussion from the previous section, we compute the phase shift at time $T_f$, for a given $M$ and $b$, by calculating the angle between the tangent vector $\T$ at times $t=0$ and $t=T_f$. Note that the phase shift decreases as $M$ increases, and as $M\to\infty$, it converges to $2\pi b/M^2$. In Figure \ref{fig:rot-shft-all}, we have compared the corresponding phase shift with the quantity $2\pi b/M^2$, taking $b=0.4$, $M=3,4,\ldots,20$, $N/M=7680$. The left-hand side shows the phase shift values for different $M$ values, whereas, on the right-hand side, we have computed the relative and absolute errors. The results also suggest that, as $M$ grows larger, the amount of both the Galilean and the phase shifts decreases at the end of one time-period.

\subsection{Case with $b \to 1^-$}\label{b=1}
 
From (\ref{eq:theta_expression}), it follows that, as $\theta_0 \to 2\pi/M$, $b\to 1^{-}$, but, in Section \ref{sec:prob-def}, we mentioned that we can have both space and time periodicity only when $\gamma=1$ (or $\theta_0 =2\pi/M$). Since, numerically, $\theta_0$ cannot be exactly $2\pi/M$ (or $b=1$), we have taken $b=1-10^{-5}$, and observed the evolution for $M^2$ time periods, i.e., until $t=2\pi$. On the other hand, as $b$ approaches $1$, the speed of the center of mass $c_M$ tends to $0$, and this implies that the vertical movement (or the third component $X_3(0,t)$, which, after removing the vertical height, is $2\pi/M$-periodic), is very small compared to the other two components. Hence, in order to understand the behavior of $\X(0,t)$ in the XY-plane, we consider the following stereographic projection in $\mathbb{C}$:
\begin{eqnarray}\label{eq:z_t-b1}
z(t) = \frac{ X_1(0,t)}{1+ \tilde X_{3}(t)} + i \frac{X_2(0,t)}{1+ \tilde X_{3}(t)}, \quad t \in [0,2\pi],
\end{eqnarray}

\noindent which is almost $2\pi$-periodic. As done for $b\in(0,1)$, we approximate the Fourier expansion of $z(t)$, and note that the dominating points in the fingerprint correspond to the squares of those integers belong to the set
\begin{equation}
\label{e:AM}
A_M = \{1\}\cup\{n\, M \pm 1\,|\, n\in\mathbb N\}.
\end{equation}

\noindent This motivates us to compare $z(t)$ and
\begin{equation} \label{eq:phi_M}
\phi_M(t) =  \sum\limits_{k \in A_M} \frac{e^{2\pi i k^2 t}}{k^2}, \quad t \in [0,1].
\end{equation}

\begin{figure}[!htbp]\centering
	\includegraphics[width=0.496\textwidth, clip=true, align=c]{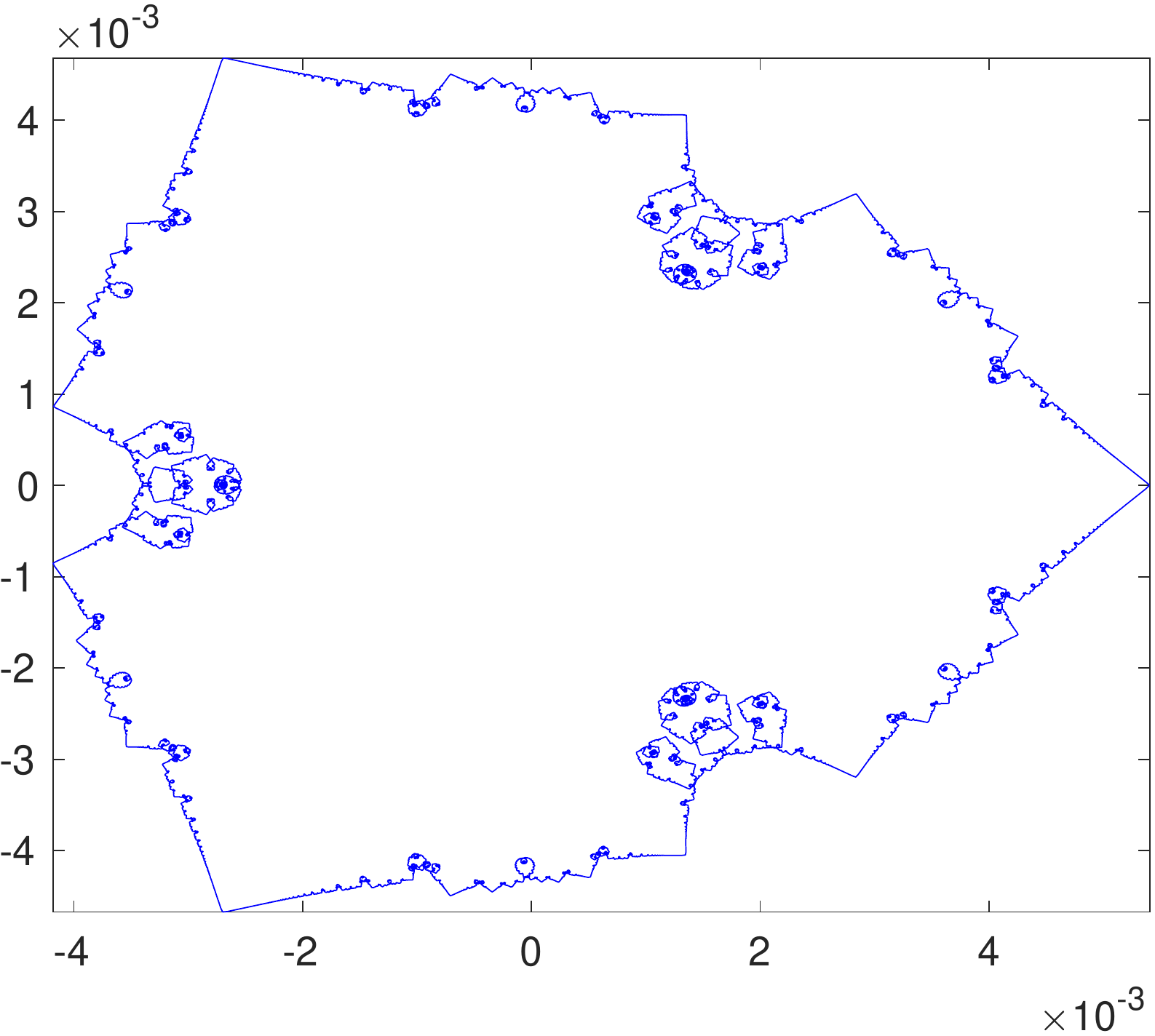}\includegraphics[width=0.504\textwidth, clip=true, align=c]{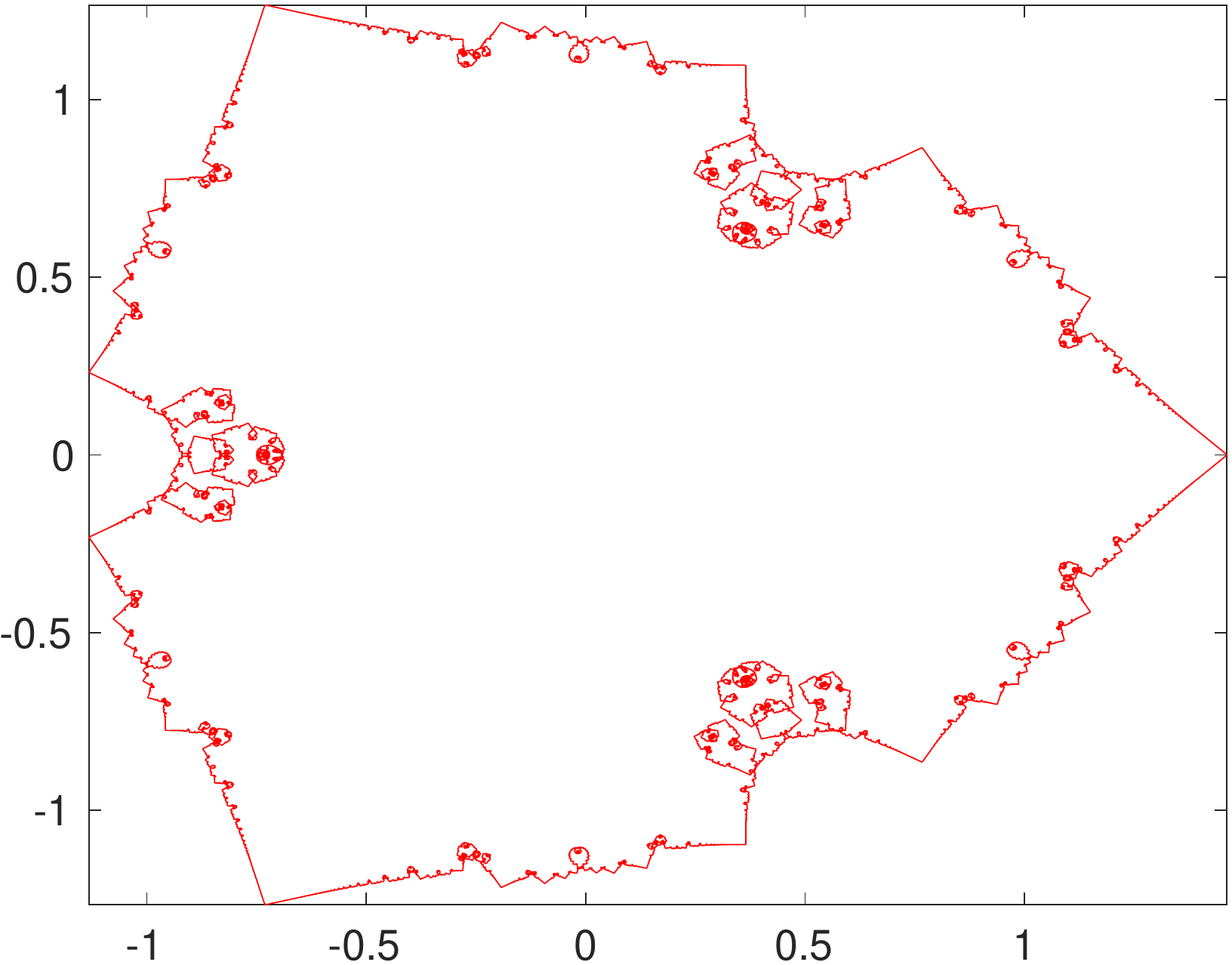}
	\caption{Left: Trajectory of a single point $z_M(t)$, for $M=3$, $N/M=2^{10}$. Right: $\phi_M(t)$ in (\ref{eq:phi_M}), where $\dim(A_M)=2^{10}$. Both curves have been computed at $N_t+1$ points in their respective domain intervals. }
	\label{fig:X0-phiM-3}
\end{figure}

\noindent In order to determine the correct orientation of the almost closed curve $z(t)$, we rotate it clockwise by an angle of $\pi/2-\pi/M$ radians, and call the resulting curve $z_M(t)$. Figure \ref{fig:X0-phiM-3} shows $z_M(t)$ (blue) and $\phi_M(t)$ (red), for $M=3$. Except for a scaling, both curves are visually the same. 

\begin{figure}[!htbp]\centering
	\includegraphics[width=0.505\textwidth, clip=true]{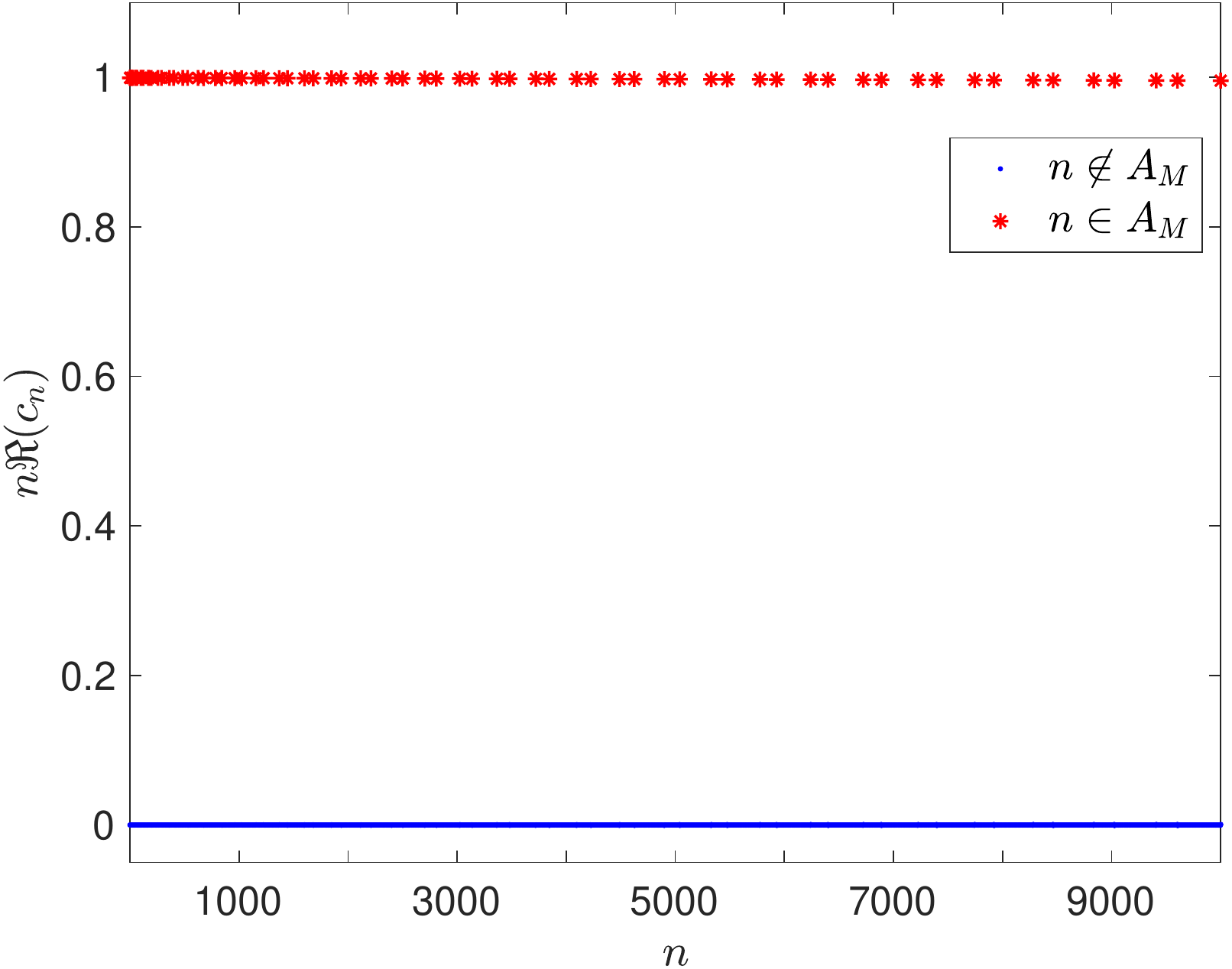}\includegraphics[width=0.495\textwidth, clip=true]{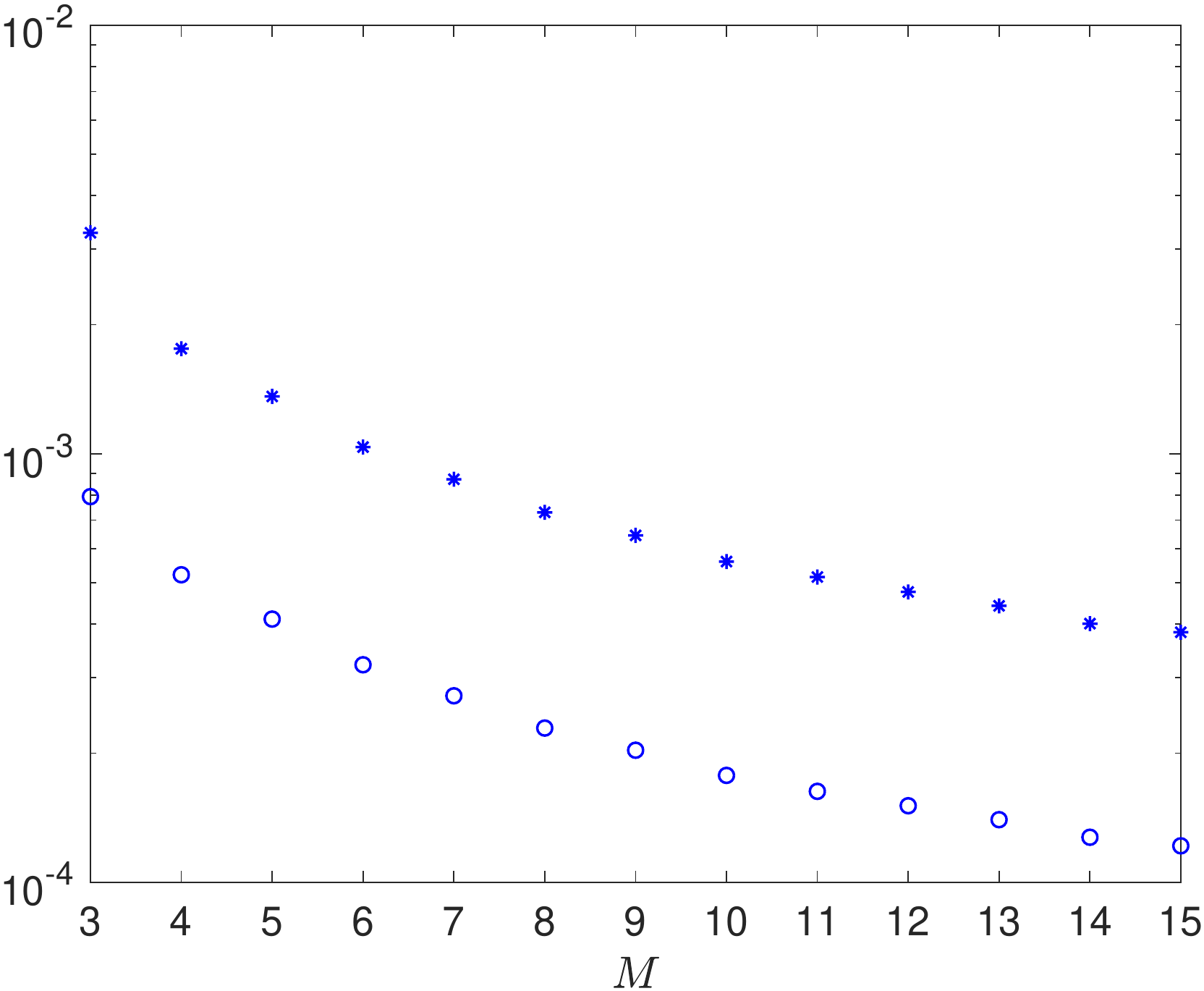}
	\caption{Left: Plot of  $n\Re{(c_n)}$ against $n$, for $M=3,N/M=2^{10}$, $b=1-10^{-5}$, $t \in [0, 2\pi]$. The dominating points (red starred) are of the form $k^2 c_{k^2}$, where $k$ is such that $\bmod(k\pm1,M)=0$; we have taken $k\in \{1,2,4,5,7,\ldots,34\}$. Right: errors $\max_t|(\phi_M - \lambda_M  z_M - \mu_M)|$ (circled) and $\max_t|(\phi_M - \lambda_M  z_M - \mu_M)/ \phi_M|$ (starred), where $\lambda_M$ and $\mu_M$ are computed from $(73)$ in \cite{HozVega2014}.}
	\label{fig:FPb1-X0-phiM-errors}
\end{figure}

In order to strengthen our claim, we take different values of $M$, $M=3,4,\ldots,15$, and compare them with the corresponding $\phi_M(t)$. In particular, we compute $\phi_M(t) - \lambda_M z_M(t) - \mu_M$, for some $\lambda_M \in \mathbb{R}$ and $\mu_M \in \mathbb{C}$ obtained using a least-square fitting, as in \cite[(73)]{HozVega2014}. Figure \ref{fig:FPb1-X0-phiM-errors} shows the absolute error ($\max_t|(\phi_M - \lambda_M  z_M - \mu_M)|$) and relative error ($\max_t|(\phi_M - \lambda_M  z_M - \mu_M)/ \phi_M|$) between $z_M$ and $\phi_M$. Note that, for different values of $M$, the time period $t\in[0,2\pi]$ is the same, so the length of the vector $z_M(t)$ would be $M^2N_t+1$, and would vary with $M$. Therefore, in order to keep a fair comparison among all the $M$ values, we have kept $N/M$ constant, i.e., $N/M=2^{10}$. As $M$ increases, the length of $z_M(t)$ increases, so we restrict ourselves to the case with $M=15$, where the length of $z_M(t)$ is $1.3608\cdot10^{8}+1$. Moreover, the plot clearly shows that the convergence is quite strong in the sense that, as $M$ increases, the discrepancy between $z_M$ and $\phi_M$ decreases. 

Finally, in Figure \ref{fig:FPb1-X0-phiM-errors}, we plot the fingerprint of the scaled $z_M(t)$, for $M=3$, $b = 1-10^{-5}$, $ t \in[0,2\pi]$. Observe that, as $b$ tends to $1$, the dominating points in the fingerprint approach $1$; or in other words, we conjecture that
\begin{equation*}
\lim\limits_{b\to 1^{-}} |n \, c_{n}|  = 
\begin{cases}
1, & \text{if } n \in A_{M}, \\
0, & \text{otherwise},
\end{cases} 
\end{equation*}

\noindent where $c_n$ are the Fourier coefficients of $z_M(t)$, which implies that the curve $z_M(t)$ converges to $\phi_M(t)$, as $b\to1^-$. 

\section{$\T(s,t_{pq})$, for $q\gg1$}
\label{sec:Talg-irrational}

In the one-corner problem with fixed boundary conditions \cite{delahoz2007}, and in the planar $M$-polygon case with periodic boundary conditions \cite{HozVega2014}, a fractal-like phenomenon was observed in the tangent vectors, too. This suggests looking for a similar behavior in the case of helical $M$-polygons. We use the algebraically constructed $\Talg$, which is correct except for a rotation about the $z$-axis, does not exhibit the Gibbs phenomenon, and can be obtained without numerical simulations.

\begin{figure}[!htbp]\centering
	\includegraphics[width=0.561\textwidth, clip=true]{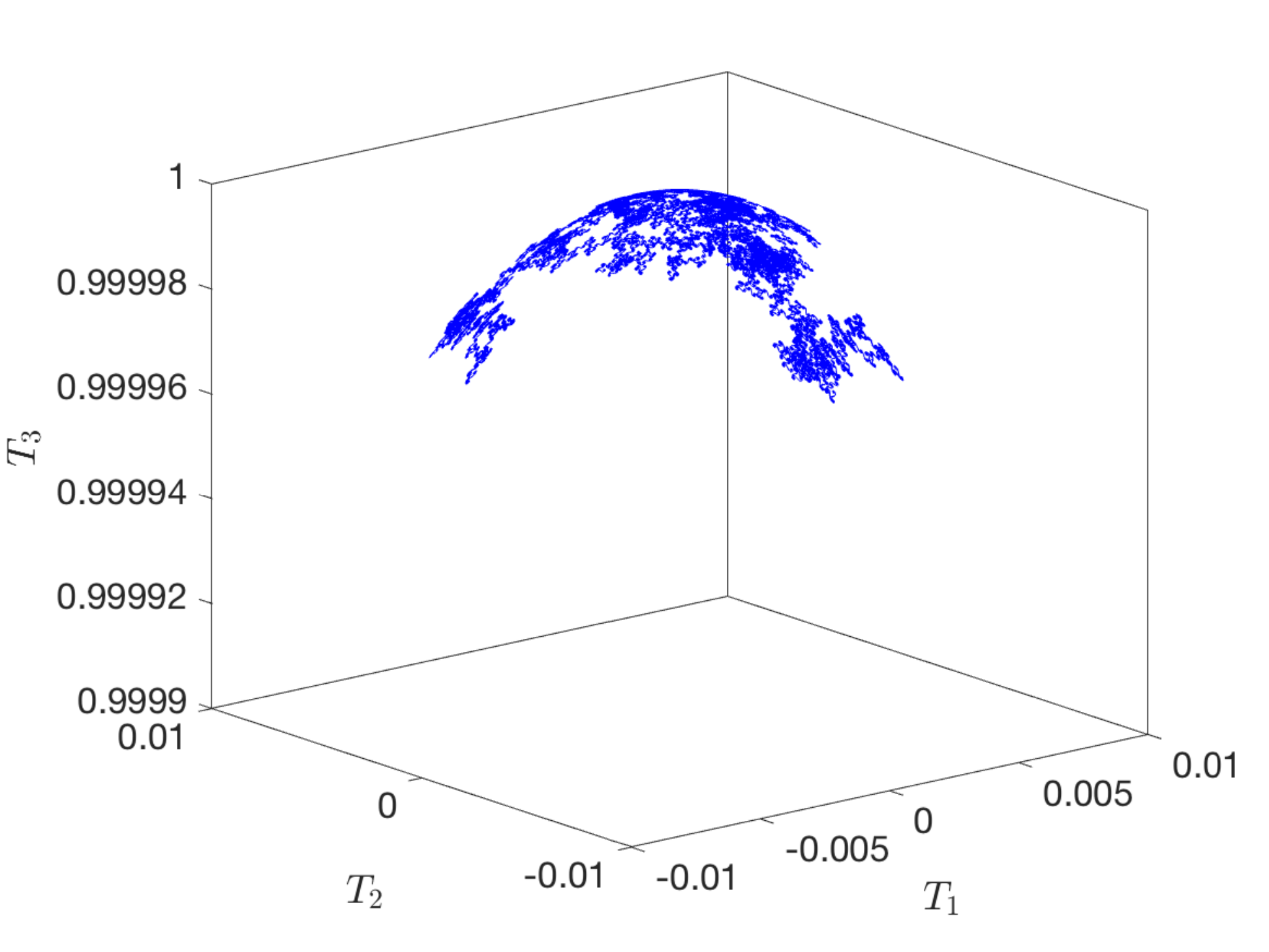}\includegraphics[width=0.439\textwidth, clip=true]{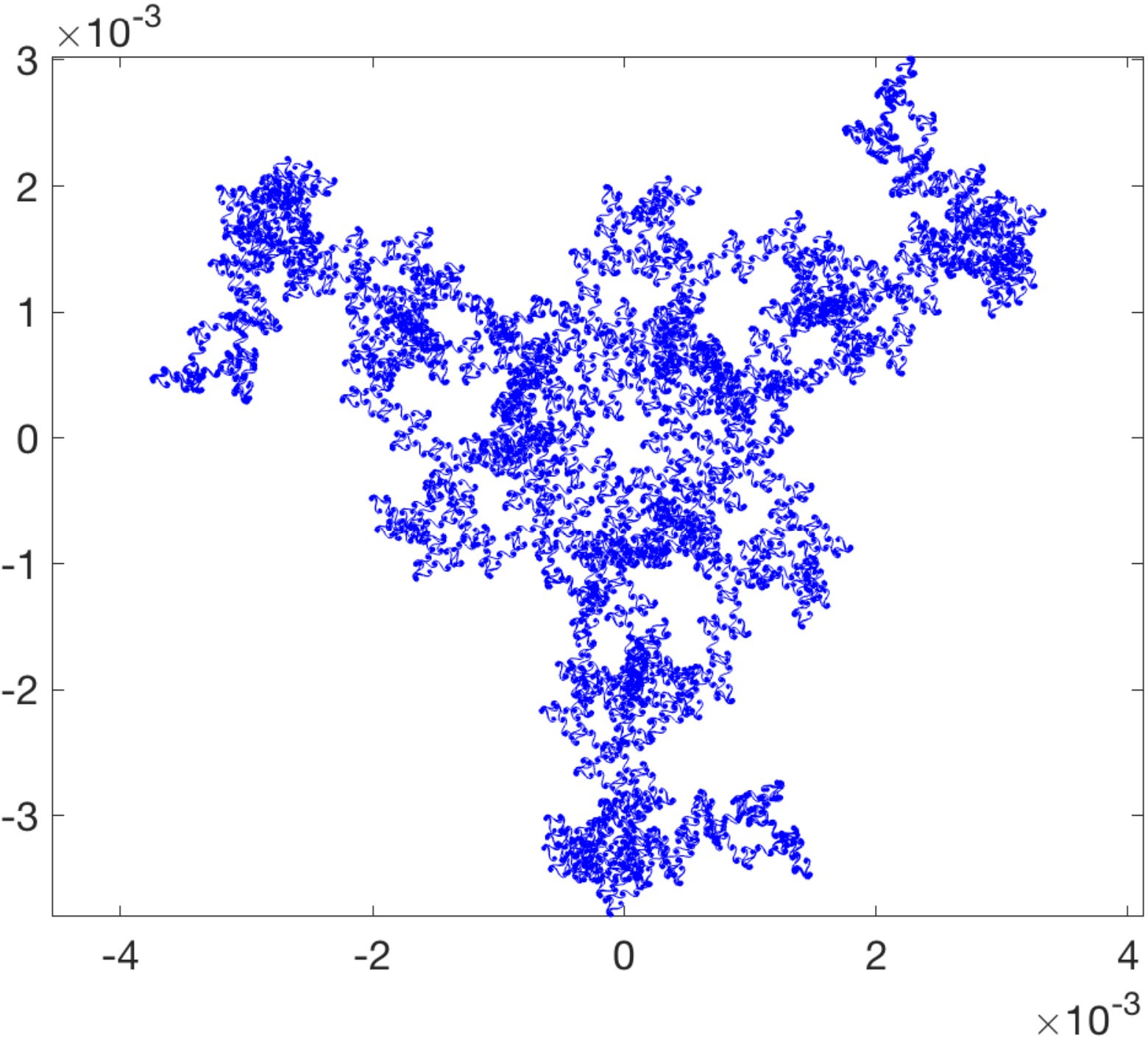}
	\caption{Left: $\Talg$, for $M=3$, $b = 1-10^{-5}$, at $t_{pq} = \frac{2\pi }{9} (\frac{1}{4} + \frac{1}{41} + \frac{1}{401})= \frac{2\pi }{9} \cdot \frac{18209}{65764}$. As $b$ tends to $1$, the values of the tangent vector concentrate around the north pole of $\mathbb{S}^2$. Right: The corresponding stereographic projection, which converges to a triskelion with the same scale spirals.} 
\label{fig:M3_IRR-b1e4}
\end{figure}

As in \cite{HozVega2014}, we take rational times $t_{pq}$, such that $q$ is very large, and there is no pair $(\tilde{p},\tilde{q})$ where both $\tilde{q}$ and $|p/q - \tilde p /\tilde q|$ are small. In particular, we take $M=3$, $t_{pq} = \frac{2\pi }{9} (\frac{1}{4} + \frac{1}{41} + \frac{1}{401})= \frac{2\pi }{9} \cdot \frac{18209}{65764}$, and note that, as $b$ moves from 0 to 1, the $Mq/2=98646$ values of $\T$ tend to concentrate on the upper half of the sphere, whereas, when $b\approx 1$, they lie very close to its north pole (see the left-hand side of Figure \ref{fig:M3_IRR-b1e4}). On the other hand, the stereographic projection seems to be even more interesting: for $b=0$, the fractal spiral-like structures appear to be at three or four different scales (see \cite[Fig. 8]{HozVega2014}), but, as $b$ approaches $1$, the equally complex structures form a shape resembling a triskelion (see the right-hand side of Figure \ref{fig:M3_IRR-b1e4}).

\section{Numerical relationship between the $M$-corner problem with nonzero torsion and the one-corner problem}\label{sec:M-1-relation}

Following the steps in \cite{HozVega2018}, we claim that the $M$-corner problem with nonzero torsion can be explained as a superposition of $M$ one-corner problems for infinitesimal times. Thus, in order to compare both cases when $\theta_0>0$, we obtain the orthonormal basis vectors, denoted as $\T_{c_\theta}$, $\nn_{c_\theta}$, $\bb_{c_\theta}$, and the curve $\X_{c_\theta}$, at $t=t_{1,q}$, $q\gg1$, by integrating \cite[(6)-(7)]{HozVega2018}. We take $c_{\theta,0}$ as defined in \eqref{e:ctheta0}, so the inner angle between the asymptotes of the tangent vectors, i.e., $\lim_{s\to-\infty}  \T_{c_\theta}(s) = \A^{-}$, $\lim_{s\to\infty}  \T_{c_\theta}(s) = \A^{+}$ (where $\A^{\pm} = (A_1,\pm A_2, \pm A_3)$ \cite{GutierrezRivasVega2003}), is equal to the angle between any two adjacent sides of the helical $M$-polygon. 

At this point, we have to rotate $\X_{c_\theta}$, $\T_{c_\theta}$, $\nn_{c_\theta}$ and $\bb_{c_\theta}$, in such a way that the rotated vectors $\X_{rot}$ and $\T_{rot}$ match the $M$-corner problem, where $\X_{rot}\equiv \M \cdot \X$ and  $\T_{rot}\equiv \M \cdot \T$, for a certain rotation matrix $\M$, which is determined by imposing that
$\T_{rot}^- = \lim_{s\to-\infty}  \T_{rot}(s) = (a\cos(2\pi / M ), -a \sin(2\pi / M), b)^T$,
$\T_{rot}^+ =\lim_{s\to+\infty}  \T_{rot}(s) = (a,0, b)^T$,  with $a^2+b^2=1$. One way of obtaining $\M$ is by writing $\M = \M_2 \cdot \M_1$, where $\M_1$ performs a rotation of an angle $\arccos(\A^+ \cdot \T_{rot}^+)$ about an axis orthogonal to these two vectors. Denoting $\tilde{\A}^- = \mathbf{M}_1 \cdot \A^-$, $\tilde{\T}^-_{rot} = \mathbf{M}_1 \cdot \T_{rot}^-$, $\M_2$ is a rotation about the axis $\T_{rot}^+$ of an angle between $\tilde{\A}^-$ and $\tilde{\T}^-_{rot}$. Bearing in mind this, $\X_{rot} = \X_0 + \M \cdot \X_{c_\theta}$, where $\X_0= (-a\pi/M,-a\pi/(M\tan(\pi/M),0)^T$ is the corner of the helical $M$-polygon corresponding to $s = 0$, at $t=0$. 

\begin{figure}[!htbp]\centering
	\includegraphics[width=0.338\textwidth, clip=true, align=t]{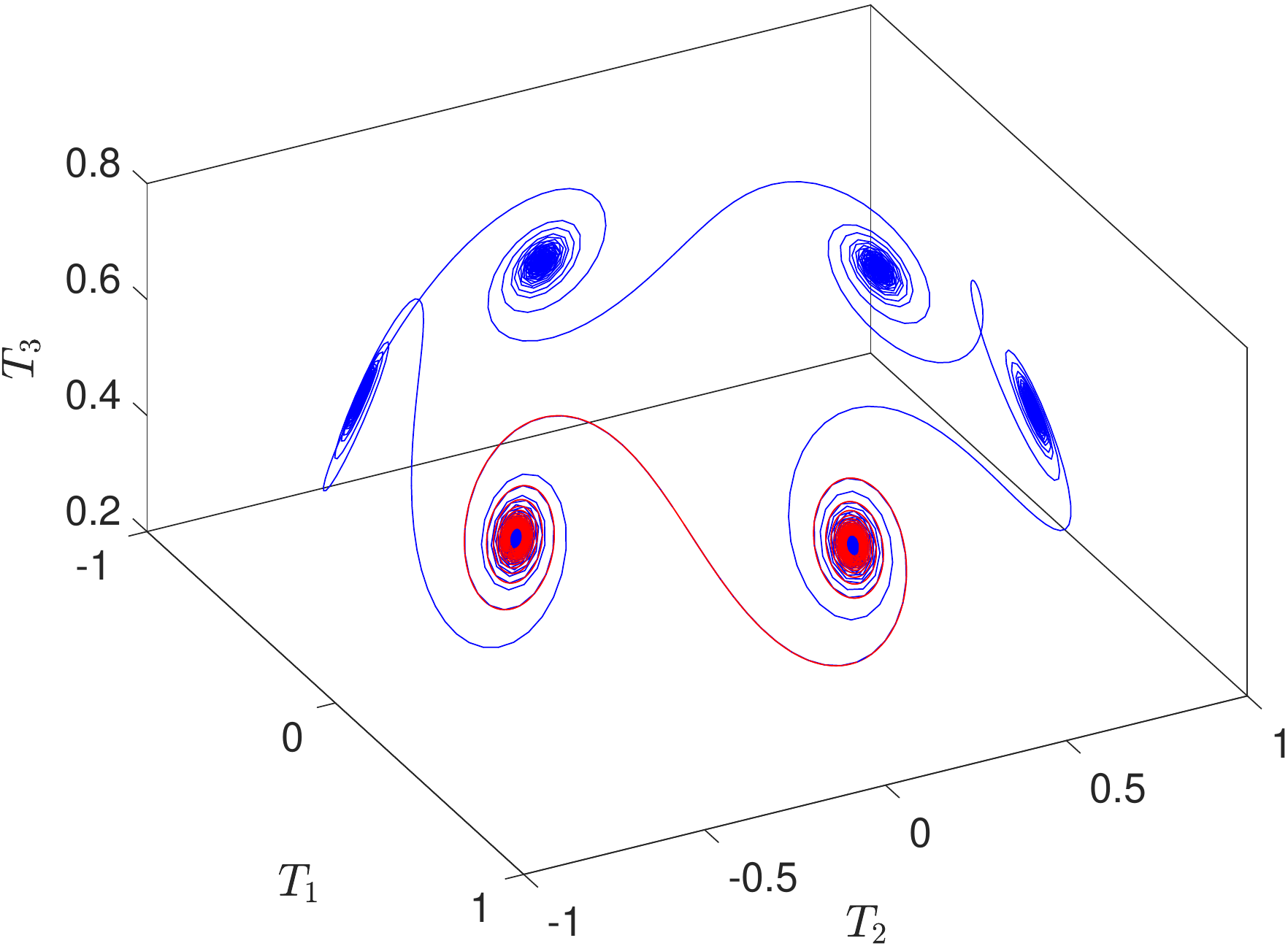}\includegraphics[width=0.32\textwidth, clip=true, align=t]{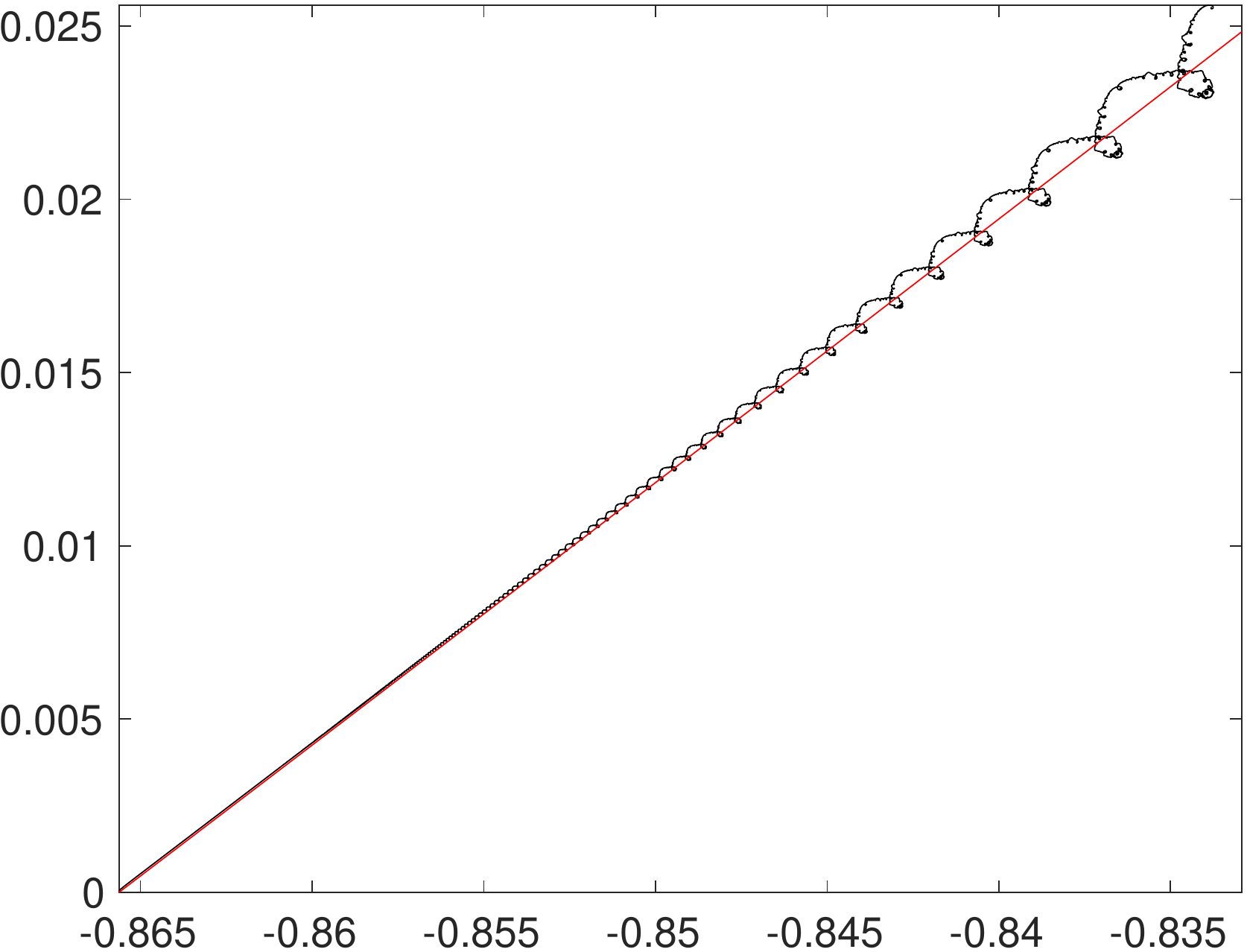}\includegraphics[width=0.342\textwidth, clip=true, align=t]{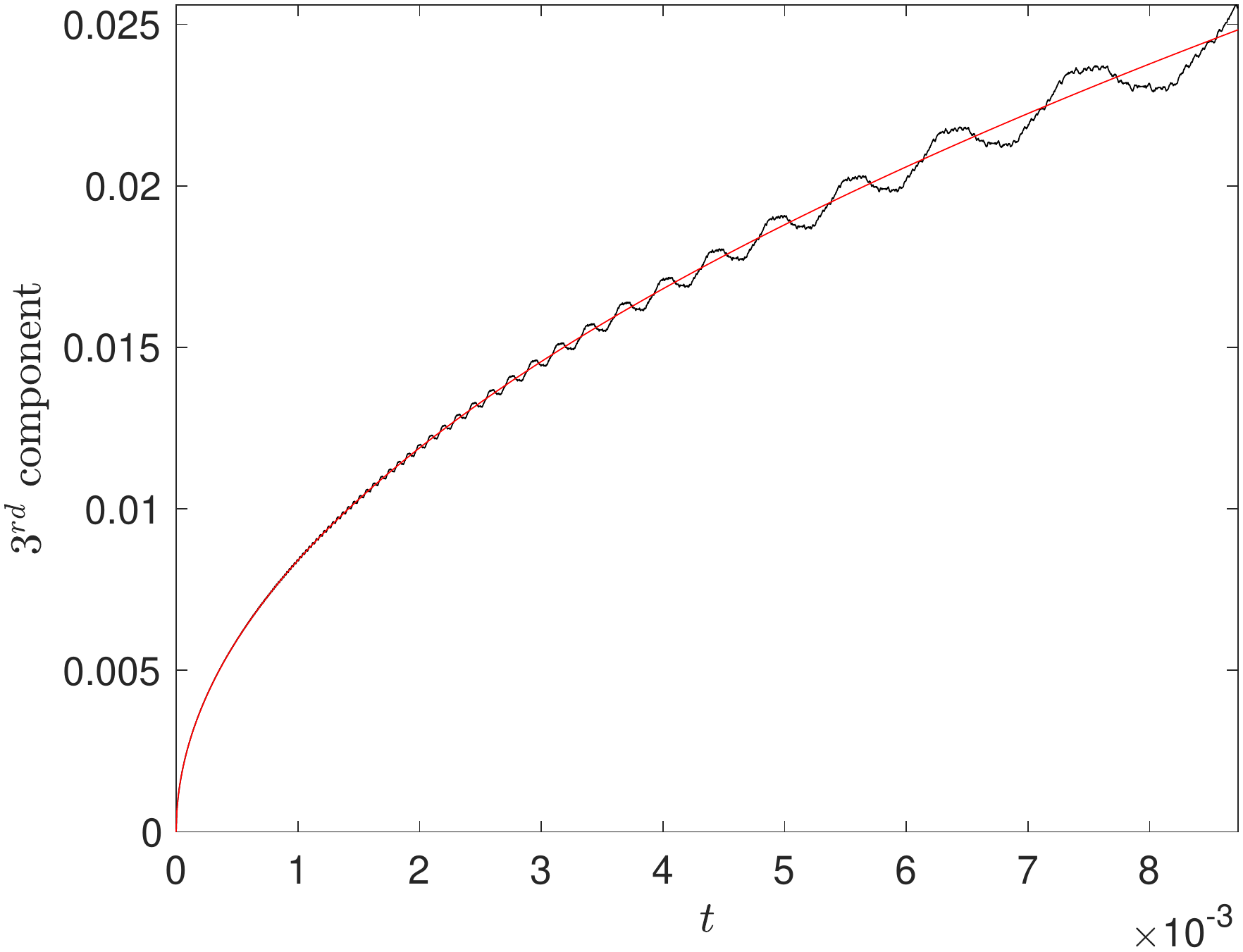}
	\caption{Left: $\T$ for the $M$-polygon problem (blue) and $\T_{rot}$ (red), for $M=6$, $\theta_0 = \pi/5$, at $t=t_{1,502}$. Center: Comparison between $\X(0,t)$ (black) and $\X_{rot}(0,t)$ (red). Right: Third component $X_{3}(0,t)$ (black) and $X_{rot,3}(0,t)$ (red), for $t\in[0,t_{1,20}]$.}	
	\label{Fig:M6_t120}
\end{figure}

In our numerical simulations, we take $\Delta s=\pi/M^2q$, and integrate \cite[(6)-(7)]{HozVega2018} at $t = t_{1,q}$, with a fourth-order Runge-Kutta method. Then, using the same $\Delta s$, we compute the $\T(s,t_{1,q})$ corresponding to the $M$-polygon problem, for $M=6$, $\theta_0 = \pi/5$, $q=502$. On the other hand, we compare the evolutions of $\X_{rot}(0,t)$ and $\X(0,t)$, for $t\in [0,t_{1,20}]$, and note that, when projected on the complex plane, $\X(0,t)$ can be very well approximated by $\X_{rot}(0,t)$, for small values of $t$. Moreover, there is also a similarity between the third components of both curves; these observations are shown in Figure \ref{Fig:M6_t120}.

\subsection{Approximation of the curvature at the origin}

From \cite{GutierrezRivasVega2003}, the curvature at $s=0$ and $t>0$ is given by $c_0(t) = \sqrt{t} | \Ts(0,t)|$. Hence, in the case of regular $M$-polygons, it can be written as $c_{\theta,0} = \sqrt{\tpq} | \Ts(0,\tpq)|$ at rational times $\tpq$. We approximate the first derivative using finite differences, as in \cite{HozVega2018}. Without loss of generality, we assume that $p=1$ and $q\equiv2\bmod 4$, so the tangent vectors are continuous at $s=0$ and $s=\pm \Delta s$, with $\Delta s = 4\pi/Mq$. Moreover, since $\theta_0 < 2\pi/M$, the Galilean shift satisfies $s_{1q} < 2\pi/Mq$, which implies that $\Talg(s_{1q})=\Talg(0)$ and $\Talg(4\pi/Mq+s_{1q})=\Talg(4\pi/Mq)$. Therefore, we approximate $c_{\theta,0}$ as
\begin{equation}\label{eq:c0_t_approximation}
c_{\theta,0} = \lim_{\substack{q \to \infty \\ q\equiv 2 \bmod 4}} \sqrt{t_{1q}} \frac{|\Talg( \frac{4\pi}{Mq},t_{1q}) - \Talg( -\frac{4\pi}{Mq},t_{1q})|}{2\cdot \frac{4\pi}{Mq}},
\end{equation}
and, after making $q \to \infty$, we recover \eqref{e:ctheta0} (see \cite[Section 4]{HozVega2018} for the intermediate steps). Then, after computing $c_{\theta,0}$ analytically from (\ref{eq:c0_t_approximation}), we have also approximated its value numerically, taking $M=6$, $\theta_0=\pi/5$, and  $q=1002,2002,\ldots,128002$. Table~\ref{table:c_0t_error} shows the discrepancies between the algebraic and numerical values. The results show clearly that, when roughly doubling $q$, the errors are approximately halved, suggesting a convergence order of $\mathcal{O}(1/q) = \mathcal{O}(t_{1,q})$.

\begin{table}[!htbp]
	\centering
	\begin{tabular}{|c|c||c|c|}
		\hline
		$q$ & $|c_{\theta,0} -  \text{approx}(c_{\theta,0})|$ & $q$ & $|c_{\theta,0} -  \text{approx}(c_{\theta,0})|$
		\\
		\hline
		1002 & $6.8511\cdot 10^{-5}$ & 16002 & $4.2878\cdot 10^{-6}$
			\\
		2002 & $3.4280\cdot 10^{-5}$ & 32002 & $2.1443\cdot 10^{-6}$
			\\
		4002 & $1.7146\cdot 10^{-5}$ & 64002 & $1.0720\cdot 10^{-6}$
			\\
		8002 & $8.5747\cdot 10^{-6}$ &	128002 & $5.3681\cdot 10^{-7}$\\
		\hline  
	\end{tabular}	
	\caption{$|c_{\theta,0}  -|\Talg(\Delta s,t_{1q})-\Talg(-\Delta s,t_{1q})|/(2\Delta s)|$, for $M=6$, $\theta_0=\pi/5$. The errors decrease as $\mathcal{O}(1/q) = \mathcal{O}(t_{1,q}).$} 
	\label{table:c_0t_error}
\end{table}

\section*{Acknowledgments}

This paper was partially supported by an ERCEA Advanced Grant 2014 669689 - HADE, by the MICINN projects MTM2014-53850-P, PGC2018-094522-B-I00 and SEV-2013-0323, by the Basque Government research group grants IT641-13 and IT1247-19, and by the Basque Goverment BERC program.

\end{document}